\documentclass[11pt]{amsart}
\usepackage{mathabx}
\usepackage{graphicx}
\usepackage{subfig}
\usepackage{tabularx}
\usepackage{wasysym,mathtools}
\usepackage{algorithm} 
\usepackage{algorithmicx}

\usepackage{lscape} 

\usepackage{mathrsfs}     
\usepackage{helvet}         
\usepackage{courier}        
\usepackage{type1cm}      
\usepackage{color,colortbl}
\usepackage{url}

\usepackage{geometry,calc,color}
\usepackage{amsmath,amssymb}

\usepackage{enumerate}
\usepackage{arydshln}
\usepackage{siunitx}
\usepackage{booktabs}
\usepackage{multirow}
\usepackage{tikz}
\usepackage{pgfplots}

\usepackage{capt-of}

\pgfplotsset{compat=1.9}

\usepackage{textcomp}

\tikzset{immagine/.style={%
		above right, inner sep=0pt, outer sep=0pt}}
    
\newcommand{\roundPrecision}{2}

\newtheorem{theorem}{Theorem}[section]

\theoremstyle{definition}

\newtheorem{problem}[theorem]{Problem}
\newtheorem{assumption}[theorem]{Assumption}
\newtheorem{comment}[theorem]{Comment}
\theoremstyle{remark}

\numberwithin{equation}{section}

\newcommand{\uh}{u_h}
\newcommand{\vh}{v_h}
\newcommand{\wh}{w_h}

\newcommand{\VKh}{V_h^K}

\newcommand{\nref}[1]{(\ref{#1})}

\newcommand{\Th}{\mathcal{T}_h}

\newcommand{\EK}{\mathcal{E}^K}

\renewcommand{\S}{\mathcal{S}}

\newcommand{\Svem}{S_a^K}

\newcommand{\Pinabla}{\Pi^\nabla_K}

\pgfplotscreateplotcyclelist{mycolorlist}{%
	blue,every mark/.append style={fill=blue!80!black},mark=*\\%
	red,every mark/.append style={fill=red!80!black},mark=square*\\%
	green!60!black,every mark/.append style={fill=green!80!black},mark=otimes*\\%
	red,mark=star\\%
	blue,every mark/.append style={fill=blue!80!black},mark=diamond*\\%
	red,densely dashed,every mark/.append style={solid,fill=red!80!black},mark=*\\%
	green!60!black,densely dashed,every mark/.append style={
		solid,fill=green!80!black},mark=square*\\%
	black,densely dashed,every mark/.append style={solid,fill=gray},mark=otimes*\\%
	blue,densely dashed,mark=star,every mark/.append style=solid\\%
	red,densely dashed,every mark/.append style={solid,fill=red!80!black},mark=diamond*\\%
}

\begin{document}

	\title[AMG for VEM]
	{
	Algebraic Multigrid Methods for  Virtual Element Discretizations: a numerical study
	 }	
	\author[D. Prada]{Daniele Prada}
	\address{IMATI ``E. Magenes'', CNR, Pavia (Italy)}%
	\email{daniele.prada@imati.cnr.it}%

	\author[M. Pennacchio]{Micol Pennacchio}
	\address{IMATI ``E. Magenes'', CNR, Pavia (Italy)}%
	\email{micol.pennacchio@imati.cnr.it}%


%
	
	\date{\today}
	\thanks{{This paper has been realized in the framework of ERC Project CHANGE, which has received funding from the European Research Council (ERC) under the European Union’s Horizon 2020 research and innovation programme (grant agreement No 694515).   
}
}%

	\begin{abstract} 
	We investigate the performance of algebraic multigrid methods for  
	the  solution of the linear system of equations arising from a Virtual Element discretization.
	We provide numerical experiments on very general polygonal meshes for  a model
	 elliptic problem with and without highly heterogeneous diffusion coefficients and we draw 
	 conclusions regarding the efficacy of the method.
	\end{abstract}

	\maketitle
	

	\section{Introduction}
	\label{intro}

The high flexibility of polytopic grids is  
crucial  
when dealing with real life problems since it  enables the treatment of complex geometries  and simpler meshing of the domain including local mesh adaptivity and nonconforming grids without requiring any special treatment, e.g.  for hanging nodes.
This explains 
the increasing interest on numerical methods  where the discretization is based on arbitrarily shaped polytopic meshes.

Here we focus on the Virtual Element Method (VEM) which is a quite recent PDE discretization method that allows for polygonal and polyhedral meshes and can be viewed as an extension of the Finite Element Method (FEM) \cite{basicVEM, hitchVEM}. 
The name
 {\em virtual} is related to the fact that 
 operators and matrices, needed in the implementation of the method, are evaluated  by relying on an implicit knowledge of the local shape functions only, and 
are computed directly in terms of the degrees of freedom through a computable 
elementwise projector onto the space of polynomials. 

The method has already been applied and extended to a wide variety of different model problems 
\cite{beirao_parab,Antonietti_VEM_Stokes,Antonietti_VEM_Cahn,beirao_stokes,beirao_Navier_Stokes,perugia_Helmholtz,beirao_elastic,beirao_linear_elasticity,VEM_3D_elasticity,VEM_mixed,VEM_discrete_fracture,VEM_Laplace_Beltrami};  the $p$ and $hp$ versions of the 
method are discussed and analyzed in  \cite{hitchVEM, antonietti_p_VEM,beirao_hp,beirao_hp_exponential}, 
 its implementation can be found in \cite{hitchVEM} and its extension 
 to the case of curved boundaries can be found in \cite{VEM_curved_beirao,VEM_curvo}. 

In this paper we focus on the efficient solution of the linear system of equations associated with a VEM discretization of a model  elliptic problem.
 The design of computationally effective solvers  for such linear systems is a crucial phase of the overall numerical process.
So far few  works in literature have been devoted to the study 
of effective numerical algorithms for the solution of these linear systems of equations.
Attempts in the recent literature focused on the increase in the condition number of the stiffness matrix 
due either to a degradation of the quality of the tessellation and/or to the increase in the polynomial order  of the method \cite{beirao_hp,Mascotto_illcond,berrone_borio_17}.
Other works tackled the problem of the increase in the condition number resulting from refining the discretization by considering domain decomposition techniques 
  \cite{Calvo_Schwarz,Calvo_MG,FETI_VEM_2D,Prada_enumath_17,FETI_VEM_3D} or $p$ refinement multigrid methods \cite{antonietti_p_VEM}. 

We observe that designing geometric multigrid solvers for VEM for the case of $h$ refinement is not straightforward. 
In fact, for such method, even if two grids are embedded one in the other, the two corresponding spaces are not. This makes the design of coarsening and prolongation operators tricky. 
The alternative is to use the algebraic version of multigrid.
Algebraic multigrid (AMG) was introduced as a method for solving linear systems based on multigrid principles, without 
exploiting 
the problem geometry; only the connections in the matrix graph are used to determine intergrid transfer operators and to define coarse grids, see e.g. \cite{Ruge.Stuben.87,Stuben2001}. 
AMG has already proved its usefulness in various  problem types and especially those discretized on unstructured grids \cite{RobScalAMG}. 
Here we want to exploit the performance of AMG solvers when very general polygonal meshes are considered, i.e. with very general  shape  as allowed by VEM.
%
	
AMG is also commonly used as a preconditioner for a Krylov subspace method such as the conjugate gradient (CG). 
The approach of using AMG with CG is often referred to as \emph{CG accelerated multigrid (AMG/CG)} and it has the advantage that the Krylov method reduces the error in eigenmodes that are not being effectively reduced by multigrid. 
	
The aim of this paper is to investigate the effectiveness 
of Algebraic Multigrid Methods (AMG) when applied to the linear system of equations associated with the lowest  order VEM discretization.
We follow the AMG/CG approach and analyze the 
performance 
of the conjugate gradient method implemented in PETSc~\cite{petsc-web-page}  and its interfaces to different AMG preconditioners in the default setting.
%
The interest in AMG techniques is mainly related 
to their potential scalability with the size of the problem to be solved,
 in the sense that 
the number of iterations required to reach convergence for a given problem does not depend on the number of the mesh nodes. 
%
The results of the present paper show that using AMG preconditioners for the solution of the linear system of equations associated with VEM discretizations is a promising approach, in terms of both scalability and reduction of the overall computational cost.
With most of the meshes considered, AMG outperforms other classical solvers available in PETSc, in terms of computational times.
However, we also verified that, when dealing with particularly complex and challenging meshes, 
not all the AMG preconditioners considered 
preserve  scalability and those that do, 
loose most of their efficiency. 
 The use of AMG/CG 
for the solution of linear systems associated with VEM discretization based on these  meshes deserves further investigation, 
beyond the scope of this work.
%

%

The paper is organized as follows. The basic notation 
and the description of the Virtual Element Method are given in Section~\ref{sec:VEM}. 
Algebraic Multigrid Methods are briefly recalled  in Section~\ref{sec:AMG} with 
  the different codes used. 
Numerical experiments to test the performance of AMG are presented in Section~\ref{sec:numerical} and finally in Section \ref{sec:conclusion} we draw conclusions on the efficiency and robustness of the AMG preconditioners considered. 
	
%

\section{The virtual element method (VEM)}\label{sec:VEM}
Let us start by recalling the definition and the main properties of the  Virtual Element Method \cite{basicVEM}. 
To fix the ideas we focus on the following  model  elliptic problem: 
\[
-\nabla \cdot(\rho \nabla u) = f \ \text{ in } \Omega, \qquad u = 0 \ \text{ on } \partial \Omega,
\]
with $f \in L^2(\Omega)$ and $\Omega \subset \mathbb{R}^2$  polygonal domain. We assume that the coefficient $\rho$ is a scalar such that for almost all $x \in \Omega$,  $\alpha \leq \rho(x) \leq M$ for two constants $M \geq \alpha >0$. The variational formulation of such an equation reads as follows
\begin{equation}\label{variational}
\left\{
\begin{array}{l}
\text{find $u \in V:=H^1_0(\Omega)$ such that}\\[1mm]
a(u,v) = ( f , v ) \ \forall v \in V 
\end{array}\right.
\end{equation}
with
\[
a(u,v) = \int_{\Omega}\rho(x)\nabla u(x)\cdot \nabla v(x)\,dx, \qquad (f,v) = \int_{\Omega}f(x) v(x)\,dx.
\]

\

We consider a family $\{ \Th \}_h$ of tessellations of $\Omega$ into a finite number of simple polygons $K$. 
The assumptions  generally considered on each tessellation $\Th$ are the following 
 (see, e.g. \cite{basicVEM,hitchVEM}):

\begin{assumption}\label{assumption_K}
	there exist constants $\gamma_0, \gamma_1, \alpha_0, \alpha_1 > 0$ such that: 

\begin{enumerate}
	\item each element $K \in \Th$ is star-shaped with respect to a ball of radius $\geq \gamma_0 h_K$, where $h_K$ is the diameter of $K$;
	\item for each element $K\in\Th$ the distance between any two vertices of $K$ is $\geq \gamma_1 h_K$;
	\item $\Th$ is quasi-uniform, that is, for any two elements $K$ and $K'$ in $\mathcal{T}_h$ we have $\alpha_0  \leq h_K / h_{K'}\leq \alpha_1 $.
\end{enumerate}
\end{assumption}



For simplicity we assume that 
	for all $K$ there exists a constant $\rho_K$ such that the coefficient $\rho$ verifies $\rho = \rho_K$ on $K$. 
	However, we remark that the virtual element method can also be defined for the case of coefficients varying within the elements of the tessellation, see e.g.  \cite{beirao_var_coef}.

The Virtual Element discretization space is defined element by element starting from the edges of the tessellation. 
More precisely, for each  polygon $K \in \Th$ we introduce the space $\mathbb{B}_1(\partial K)$ as
\begin{gather*}
\mathbb{B}_1(\partial K) = \{ v \in C^0(\partial K): v|_{e} \in \mathbb{P}_1\ \forall e \in \EK \},
\end{gather*}
where $\mathbb{P}_1$ denotes the set of polynomials of degree less than or equal to $1$, and 
$\EK$ 
the set of edges of the polygon $K$. 
Letting \begin{equation}\label{deflocalspace}
\VKh =   \{ v \in H^1(K):\ v|_{\partial K} \in \mathbb{B}_1(\partial K),\ \Delta v =0 \text{ in } K \}
\end{equation}
the discrete space $V_h$ is then defined as
\begin{align*}
V_h &= \{ v \in V: w|_{K} \in \VKh, \forall K \in \Th \}.
\end{align*}

\

Let $(\cdot,\cdot)$ be the scalar product in $L^2$, 
$a(u,v)= (\rho \nabla u, \nabla v)$ and  $a^K$ the restriction of $a$ to $K$.
Using a Galerkin approach, we look for $u_h\in V_h$ such
that for all $v_h\in V_h$ 
$$a(u_h,v_h)=\int_\Omega f \, v_h \, dx. $$

Both terms at the right 
and at the left hand sides 
 cannot be computed exactly with the knowledge of the value of the degrees of freedom of $u_h$ and $v_h$ only.
Setting, 
for each $K \in \Th$
\[
a^K(u,v)=\int_K \rho \nabla u\cdot \nabla v 
\]
we observe that, by using Green's formula, given any $v \in V^{K}_h$ and any $p \in \mathbb{P}_1(K)$ 
\[
a^K(p,v) = - \rho_K\int_K v \Delta p + \rho_K \int_{\partial K} v \frac {\partial p}{\partial n}.
\]

Since on each edge of $K$ $v$ is a known linear and $\partial w/\partial n$ is a known constant, the right hand side can be computed exactly and directly from the degrees of freedom of $v_h$. 
This allows to define  the ``element by element'' computable projection
operator $\Pinabla : V^{K}_h \longrightarrow \mathbb{P}_1(K)$  
\[
a^K(\Pinabla  u, q ) = a^K(u,q) \quad \forall q \in \mathbb{P}_1(K), 
\]
and we clearly have 
\begin{equation}\label{1}
	a^{K}(u,v) = a^K(\Pinabla  u,\Pinabla  v) + a^K(u - \Pinabla  u, v - \Pinabla  v).
\end{equation}
The virtual element method stems 
from replacing the second term of the sum on the right hand side (which cannot be computed exactly), with an ``equivalent''  term, where the bilinear form $a^K$ is substituted by a computable symmetric bilinear form $\Svem$, resulting in defining
\[
a^{K}_h(u,v) = a^K(\Pinabla  u,\Pinabla  v) + \Svem(u - \Pinabla  u, v - \Pinabla  v) .
\]

As it is usually done in VEM, the bilinear form $a$ is then replaced with
a suitable approximate bilinear form  
$a_h : V_h \times V_h \to \mathbb{R}$ defined by
\[
a_h(u_h,v_h) = \sum_K a_h^K(u_h,v_h) .
\]

Different choices are possible for the bilinear form $\Svem$ (see \cite{beirao_stab}), the essential requirement being that it satisfies
\begin{equation}\label{2}
	c_0  a^K(v,v) \leq  \Svem(v,v) \leq c_1 a^K(v,v)\quad \forall v \in V^K_h\ \text{ with } \Pinabla  v=0,\end{equation}
for two positive constants $c_0$ and $c_1$, so that we have 
\begin{equation}\label{elemequiv}
	(1+c_0)   a^K(v,v) \leq   a^K_h(v,v) \leq (1+c_1) a^K(v,v)\quad \forall v \in V^K_h,
\end{equation}
and the local discrete bilinear forms satisfy the two following properties:
\begin{itemize}
	\item {\em stability}: 
	$$
	(1+c_0)   a^K(v,v) \leq   a^K_h(v,v) \leq (1+c_1) a^K(v,v)\quad \forall v \in V^K_h,
	$$
	\item {\em consistency}: for any $v_h\in V^K_h$ and $p\in \mathbb{P}_{1}(K)$ 
	\begin{equation}\label{mconsistency}
		a_h^K(v_h,p) = a^K(v_h,p) .
	\end{equation}
\end{itemize}
In the numerical tests performed in Section \ref{sec:numerical} we made the standard choice of defining $\Svem$ in terms of the vectors of local degrees of freedom as the properly scaled euclidean scalar product. 

As far as  the linear form $f$ on the right-hand side of the variational problem \nref{variational}, it  is discretized by
$ f_h : \VKh \rightarrow \mathbb{R}$ such that
$$f_h(\vh) :=   (\Pi^0_K f,\vh)_{0,K}\qquad  \forall K\in \Th,$$
where $\Pi^0_K : \VKh \rightarrow \mathbb{R}$ denotes the $L^2(K)$-orthogonal projection onto constants, defined for any $\wh\in \VKh$ to be such that
$$\int_K( \wh - \Pi^0_K \wh ) dx = 0.$$

Thus the virtual element discretization of (\ref{variational}) yields the following discrete problem:
\begin{problem}\label{discrete_full} Find $u_h \in V_h$ such that
	\[  a_h(u_h,v_h) = f_h(v_h) \qquad  \forall v_h \in V_h. \]
\end{problem}

\noindent For the study of the convergence, stability and robustness properties of the method we refer to~\cite{basicVEM,beirao_var_coef}.

\subsection{Matrix form}\label{sec:matrix}
We now  focus on the construction of the linear system of equations stemming from problem \ref{discrete_full}. 
Proceeding as in the Finite Element Method (FEM), on each element, we consider a system of $N$ shape functions, $\{ \varphi_i(\mathbf{x})\}_{i=1}^N$ associated with the $N$ vertices of the element. 
We assume that the shape functions satisfy the Lagrange conditions, that is, 
they take value 1 at the vertex which they are associated with and zero at the other vertices.
We deal with the standard  basis for the space $\mathbb{P}_K$ of local linear polynomials on each
element $K$, that is, the set of scaled monomials of degree 1 that are defined  on the element $K$ as follows
\begin{equation}
\mathcal{M}_K:=\left \{ m_1(x,y) := 1, \quad m_2(x,y) :=\frac{
  x - x_K }{h_K}, \quad m_3(x,y) := \frac{ y - y_K  }{h_K} 
\right \},
\end{equation}
with $x_K,y_K$ and $h_K$ 
being the barycenter and diameter of the element $K$, respectively. 

If we write the virtual element solution $\uh$ as 
$\uh = \sum_{i=1}^N dof_i(u_h)\, \varphi_i,$ 
problem \ref{discrete_full} 
can be rewritten 
 in matrix form as 
\begin{equation}\label{system}
A \mathbf{u} = \mathbf{f}, 
\end{equation}
where the elements of matrix $A$ and of $\mathbf f$ are:
\begin{align*}
a_{j,i}&=  \sum_{K\in\Th} \left(  a^K_h(\nabla \Pinabla \varphi_i,\nabla \Pinabla \varphi_j) + \S^K(\varphi_i- \Pinabla\varphi_i, \varphi_j- \Pinabla\varphi_j)\right),\quad i, j = 1,\dots, N,\\
 f_j &= \sum_K \left( \Pi^0_K f , \varphi_j \right)_{0,K},\quad j = 1,\dots, N .
\end{align*}

Both these terms are defined through sums over elements thus  
 the virtual element method implementation is similar to that of 
a standard finite element method.
The main difference with FEM is that 
the computation of the VEM local stiffness matrices relies on computing the local projector $\Pi^K$ on each element first.
Moreover, 
while the implementation of a  finite element may rely on a mapping to a reference element, this is not  possible in VEM because the mesh elements are allowed to be general  polygons. 

%
%

%

\section{Algebraic multigrid}\label{sec:AMG}

Algebraic multigrid (AMG) was introduced as a method for mimicking the performance of geometric multigrid on unstructured grids. 
The multilevel process is realized from a purely algebraic standpoint without 
any  explicit knowledge of the problem geometry. 
The most popular published implementation of the algorithm is the one by Ruge and Stuben (RS) \cite{Ruge.Stuben.87};
for an introduction to AMG, see e.g. \cite{Briggs.06,Trottenberg, Falgout.introAMG} and  \cite{Yang} for parallel AMG methods.

The  two main components of multigrid 
 are \emph{smoothing} and \emph{coarse-grid correction}. 
Coarse-grid correction involves operators that transfer information between fine and coarse ``grids'', i.e.,
from the vector space $\mathbb{R}^n$ to the coarse vector space $\mathbb{R}^{n_c}$. 
{\em Prolongation} maps the coarse grid to the fine grid and is just the 
matrix $P : \mathbb{R}^{n_c} \rightarrow \mathbb{R}^n$. {\em Restriction} 
maps the fine grid to the coarse grid and is the transpose of interpolation, $P^T$. 
The two-grid method for solving \nref{system}  is then: 

\

{\rm\tt  \begin{algorithmic}[a]
\State  do $n_1$ smoothing steps on $A\mathbf{u} = \mathbf{f}$. 
\State  Compute residual $\mathbf{r} = \mathbf{f} - A \mathbf{u} = A \mathbf{e}$. 
\State  Solve $A_c \mathbf{e}_c = P^T \mathbf{r}$. 
\State  Correct $\mathbf{u}\leftarrow \mathbf{u}+P \mathbf{e}_c$. 
\State  do $n_2$ smoothing steps on $A\mathbf{u} = \mathbf{f}$. 
\end{algorithmic}}
\

The vector $\mathbf{e} = A^{-1}\mathbf{f} - \mathbf{u}$ is the difference between the exact solution and the current iterate and it
is called the \emph{error}. 
The coarse system $A_c \mathbf{e}_c=P^T \mathbf{r}$  is solved by re-applying the algorithm, yielding a hierarchy of coarse grids, transfer operators, and coarse-grid systems. 

In the AMG method these components are realized in a purely algebraic way, without exploiting any knowledge of the underlying discretization method. The \emph{smoother} (or relaxation process)   is generally an iterative method  such as point-wise Gauss-Seidel, which effectively eliminates some errors, while other components of the error are reduced quite slowly. The smooth error component, in the case of AMG, is then defined to be any component of the error not reduced by the relaxation. The nonzero structure of the matrix is used to determine the adjacency relationships between the unknowns and the smooth error components are related  to strongly connected degrees of freedom. 
We say that $u_i$ is connected to $u_j$ if the element matrix $a_{ij}\neq0$. 
The magnitude of $a_{ij}$ indicates how strong the connection is and how much influence the error at $j$ has on the error at $i$ in the relaxation. 
Following the classic RS algorithm \cite{Ruge.Stuben.87},  $u_i$ is defined to be strongly connected to $u_j$ if 
$-a_{ij}\geq \theta \max_{k\neq i} (-a_{ik}),  $
 with $\theta$ being a user-defined strength threshold.
%

An alternative successful coarsening strategy is based on the so-called \emph{smoothed aggregation} (SA) technique \cite{Vanek_Mandel_Brezina}, where a coarse node is defined by an aggregate of a root point $i$  and all the neighboring nodes $j$  such that $a_{ij} > \theta \sqrt{|a_{ii}a_{jj}|}$.  

The AMG algorithms have setup costs associated with the automatic selection of coarse grid operators that other multilevel methods do not have in general. This automatic selection makes implementing the algorithm into a given application simpler, but it does add extra computations to every simulation that cannot be a priori estimated. 

In this work we follow the \emph{CG accelerated multigrid (AMG/CG)} approach where  AMG 
is used as a preconditioner for the conjugate gradient (CG) method. 
This approach has 
the advantage  that the Krylov method reduces the error in the eigenmodes that are not being effectively reduced by multigrid.

All our numerical experiments were carried out by using the Portable Extensible Toolkit for Scientific computation 
(PETSc)~\cite{petsc-web-page}. 
We now briefly recall the methods and the corresponding codes used for the numerical tests.

 \subsection{Methods}\label{sec:methods}

As iterative solvers, we used the conjugate gradient method implemented in PETSc~\cite{petsc-web-page} and its interfaces to different AMG preconditioners, with default settings. More precisely we considered:
\begin{itemize}
	\item GAMG: native AMG preconditioner implemented in PETSc. We tested two different versions: a classical AMG method (c-GAMG) and a smoothed aggregation AMG method (a-GAMG)~\cite{petsc-user-ref}.
	\item BoomerAMG: a parallel algebraic multigrid solver and preconditioner, which is part of the \textit{hypre} library~\cite{falgout2006} (\url{http://www.llnl.gov/CASC/hypre/}).
	\item ML: Multi Level Preconditioning Package, a smoothed aggregation algebraic preconditioner developed at Sandia National Laboratories~\cite{ml-guide} (\url{https://trilinos.org/packages/ml/}).
\end{itemize}

We also compared the performance of several direct solvers, which were called via
the interfaces available in PETSc~\cite{petsc-web-page}:
\begin{itemize}
	\item SuperLU and SuperLU\_Dist: sparse LU codes developed by Jim
	Demmel, Xiaoye S. Li, and John Gilbert~\cite{li05}
	(\url{http://crd-legacy.lbl.gov/˜xiaoye/SuperLU});
	\item UMFPACK: part of the SuiteSparse package developed by Timothy Davis~\cite{davis2004-umfpack} (\url{http://www.cise.ufl.edu/research/sparse/});
	\item MUMPS: MUltifrontal Massively Parallel sparse direct Solver, developed by
	Patrick Amestoy, Iain Duff, Jacko Koster, and Jean-Yves L’Excellent~\cite{amestoy2001,amestoy2006} (\url{http://www.enseeiht.fr/lima/apo/MUMPS/credits.html});
	\item PaStiX: Parallel LU and Cholesky solvers~\cite{henon2002} (\url{http://pastix.gforge.inria.fr/}).
\end{itemize}
SuperLU and UMFPACK are sequential solvers, whereas SuperLU\_Dist, MUMPS, and PaStiX can handle distributed memory systems. These last three solvers were run using $2$ processes.

All the experiments were run on a machine equipped with processor Intel$^\text{\textregistered}$ Core$^\text{\texttrademark}$ i7-7820HQ, operating system Ubuntu Linux 16.04 LTS, memory 64GB, 2400MHz DDR4 Non-ECC SDRAM.

\section{Numerical tests}\label{sec:numerical}
In this section we present results for  the solution of linear systems arising from the discretization with the virtual element method of degree one  of the following model problem: 
\begin{align}\label{eq:poisson}
-\nabla\cdot(\rho \nabla u) &= f\qquad \text{in }\Omega = (0,1)^2,\\
u &= 0\qquad \text{on }\partial\Omega,\nonumber
\end{align}
with $\rho$ diffusion coefficient.
Boundary conditions and loading term are chosen so that \linebreak 
$u = \frac{1}{2\pi^2}\sin(2\pi x)\sin(2\pi y)$
is the exact solution.
Several different meshes and diffusion coefficients $\rho$ are considered.
We start testing a constant-coefficient diffusion problem on a regular polygonal mesh,  then we deal with  irregular polygonal meshes, agglomeration of meshes and finally, highly heterogeneous coefficients $\rho$. 

We analyze the performance of the Conjugate Gradient method (CG) preconditioned with different  AMG preconditioners.
 All problems are run with different AMG codes but with fixed parameters. 

\begin{figure}
		\centering
			\captionof{figure}{Hexagonal mesh}
		\includegraphics[width=6.cm]{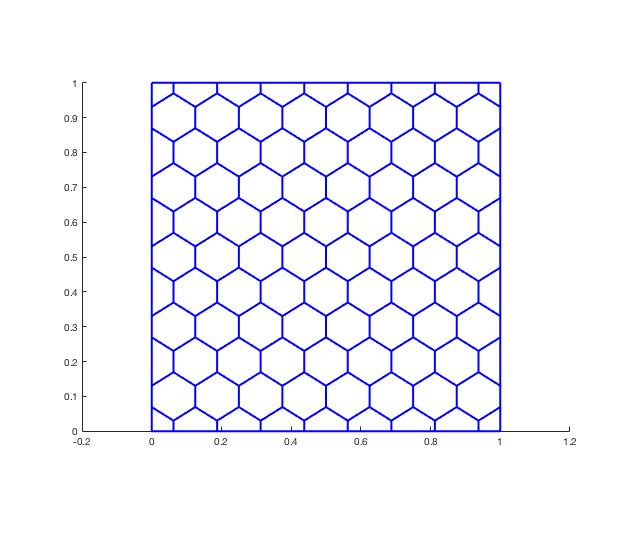}
		\label{fig:hexa}
\end{figure}

\begin{figure}
	\centering
	\caption{
	Example of Voronoi mesh. 
	}
	\subfloat{\includegraphics[width=0.33\textwidth]{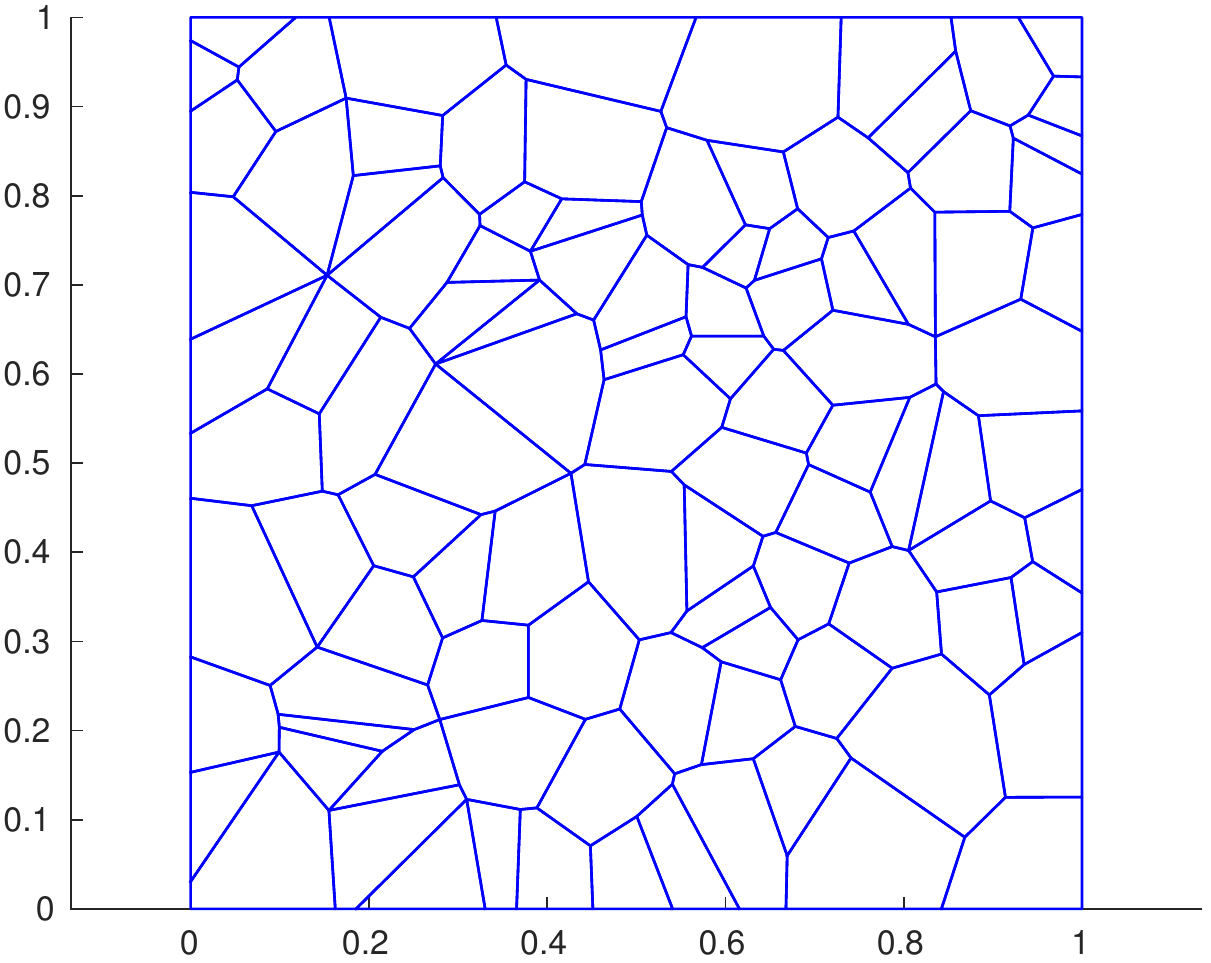}}
	\label{fig:mesh-voro} 
\end{figure}

\begin{figure}
\centering
\caption{
Example of mesh of "horse" cells. 
}
\subfloat{\includegraphics[width=0.4\textwidth]{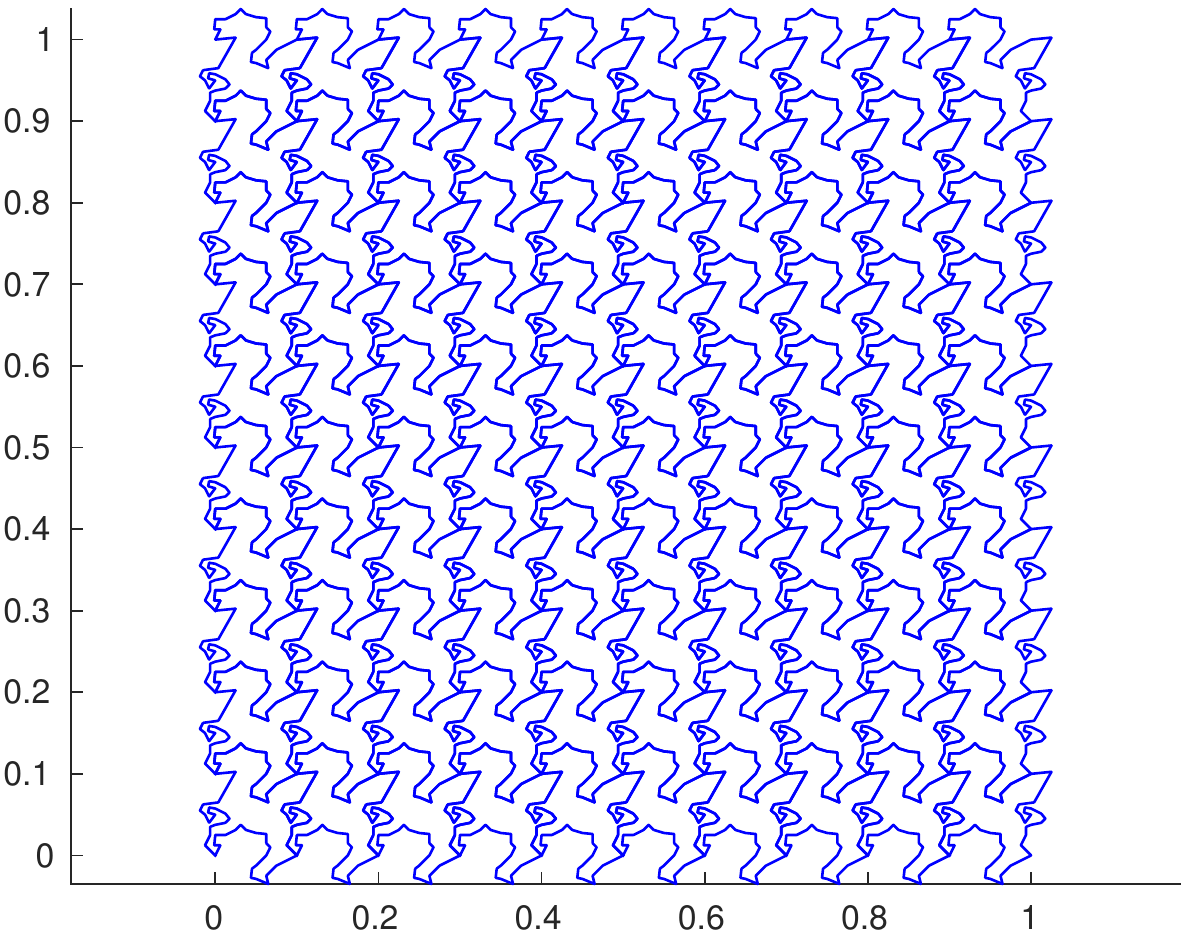}}
\label{fig:mesh-horse}
\end{figure}

\begin{figure} [h!]
	\centering
	\caption{
	Mesh obtained by embedding the third iterate of the snowflake into rectangles.}
	\begin{tikzpicture}
	\node[immagine] at (0,0) {\includegraphics[width=0.30\textwidth]{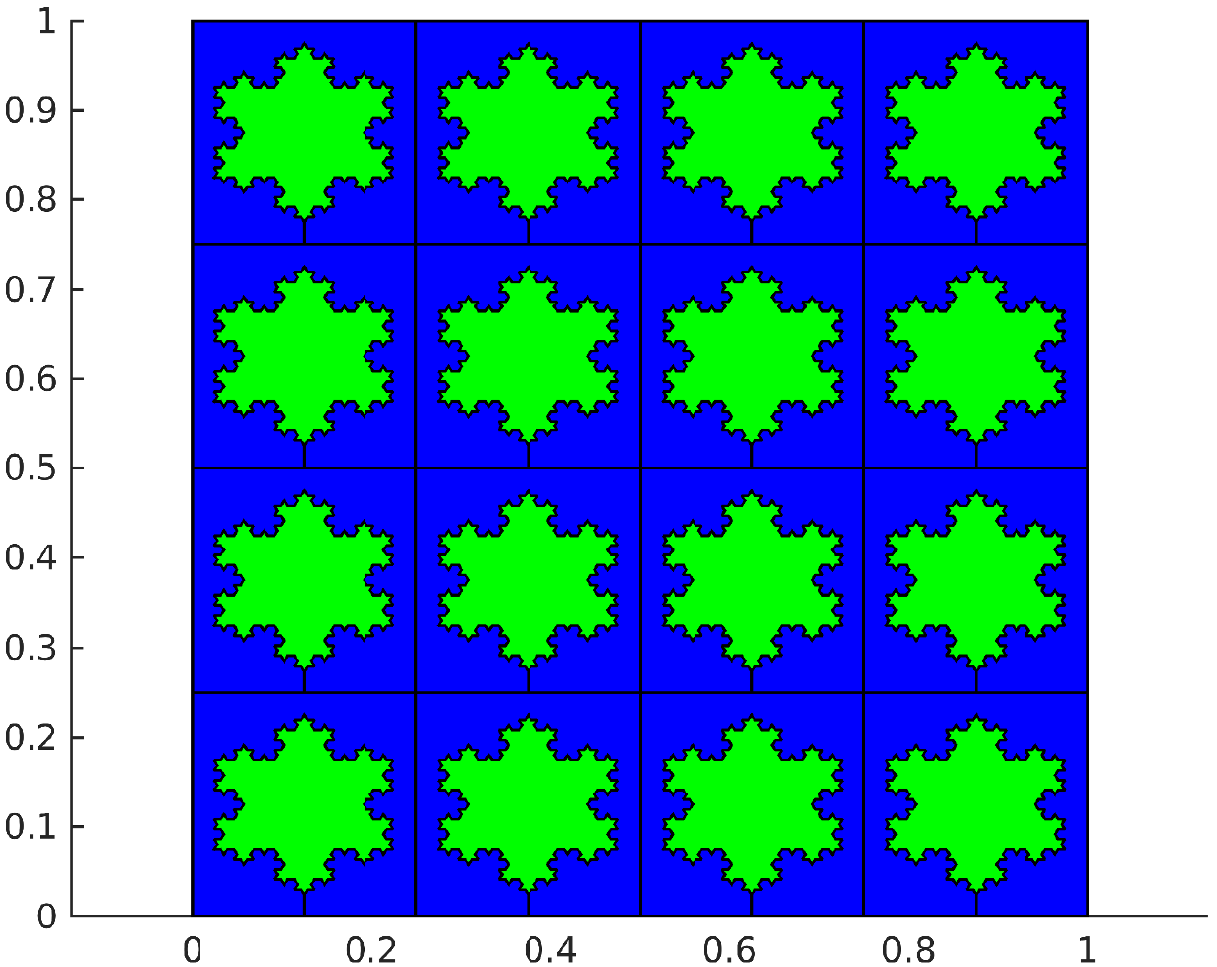}};
	\node[immagine] at (6.5,1.6) {\includegraphics[width=0.15\textwidth]{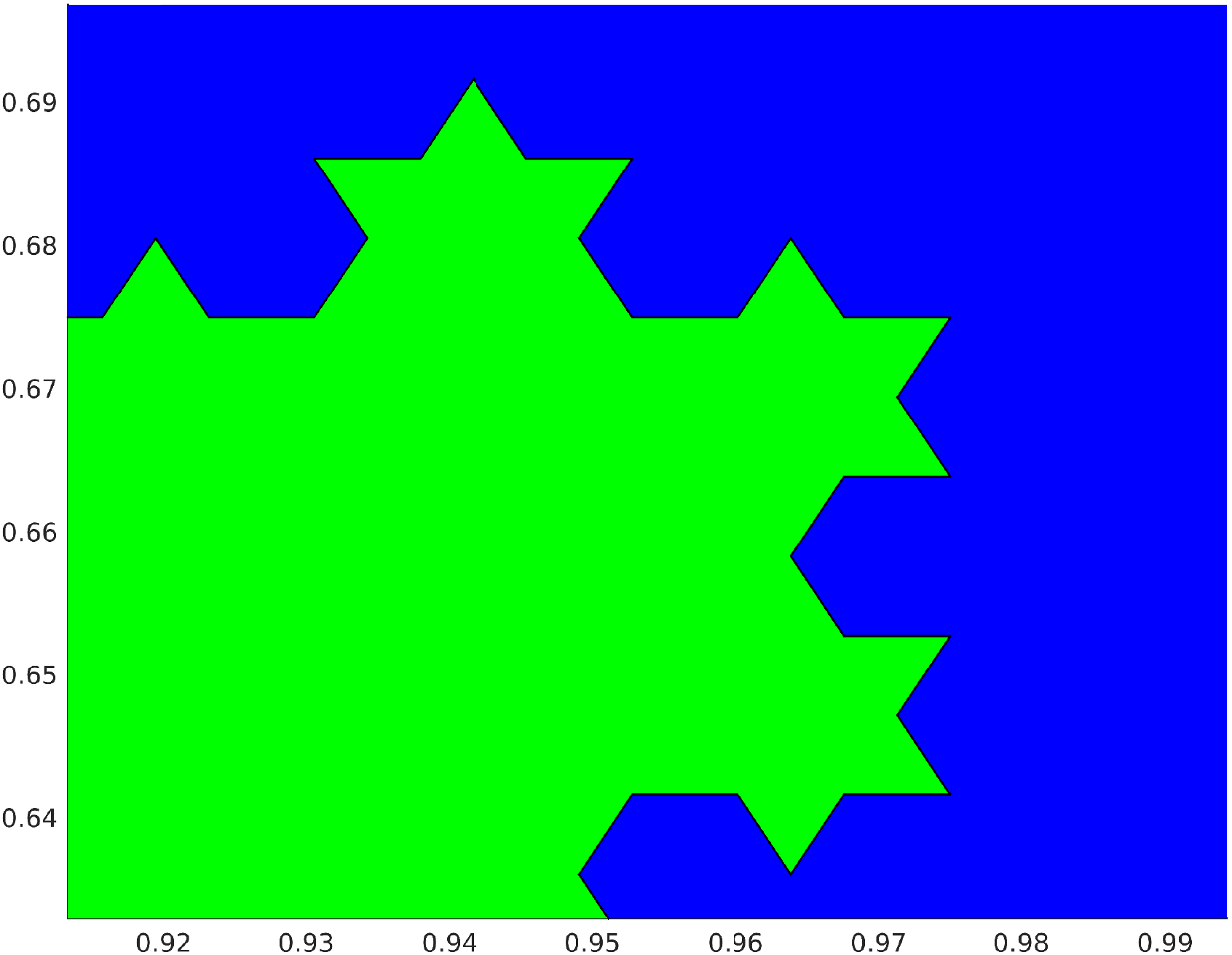}};
	\draw[ultra thick, red] (3.7,2.5) -- (4.5,2.5) -- (4.5,2.9) -- (3.7,2.9) -- (3.7,2.5);
	\draw[red] (4.5,2.5) -- (6.5,1.62);
	\draw[red] (4.5,2.9) -- (6.5,3.5);
	\end{tikzpicture}
	\label{fig:mesh-koch}
\end{figure}

 Let $\Omega$ and 
 $\Th$ be the computational domain and a polygonal tessellation, respectively. We define:
 \begin{itemize}
 	\item $N_\textup{elt}$, number of polygons of $\Th$, 
	\item $N_\textup{v}$, number of vertices of $\Th$.
 	\item $\displaystyle h = \max_{K\in\Th} h_K$, where $h_K$ is the diameter of element $K\in \Th$.
 	\item $\displaystyle h_\textup{min} = \min_{K\in\Th} h_{\textup{min},K}$, where $h_{\textup{min},K}$ is the minimum distance between any two vertices of $K$.
 	\item $\displaystyle \gamma_0 = \max_{K\in\Th}\frac{h_K}{\rho_K}$, where $\rho_K$ is the radius of the largest ball that is contained inside $K$.
 	\item $\displaystyle \gamma_1 = \max_{K\in\Th}\frac{h_K}{h_{\textup{min},K}}$.
 \end{itemize}
%
 We consider  different meshes by varying the shape of each cell and we study if the performance of the methods is affected by the shape of the cell and/or the presence of very small/large edges. 
We deal with the following different polygonal meshes $\Th$: 
\begin{enumerate}[1)]
	\item regular hexagons meshes (Figure~\ref{fig:hexa});  
	\item Voronoi meshes from uniformly random seed points 
	(Figure~\ref{fig:mesh-voro});
	\item meshes of "horse"  cell (Figure~\ref{fig:mesh-horse}); each horse is made up of 76 edges;
	\item sequence of meshes obtained by embedding successive iterates of the Koch snowflake into a rectangle  (Figure~\ref{fig:mesh-koch}).
\end{enumerate}
We remark that Voronoi and hexagonal 
meshes are  the ones 
more likely to be used, whereas all the other meshes are considered here only for stress testing of AMG preconditioners.
In particular we observe that snowflake meshes are particularly complex and challenging, as they are characterized by very small edges on the boundary of the snowflake and greater edges on the boundary of the square. These meshes may be of interest in domains with periodic structures. 

\begin{figure}
\centering
\caption{
Mesh obtained by aggregating Voronoi cells. Boundaries of aggregates are marked in red.}
\subfloat{\includegraphics[width=0.39\textwidth]{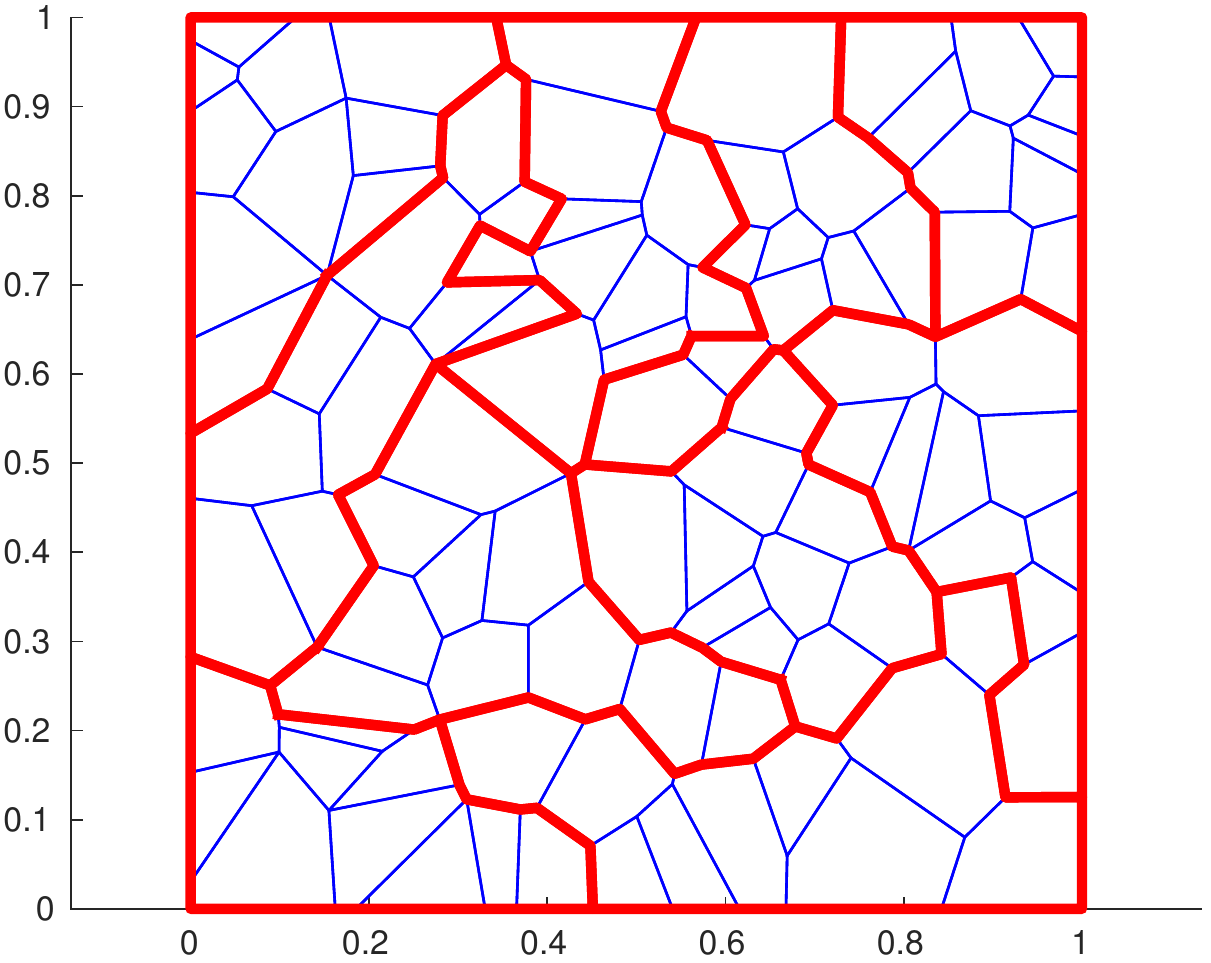}}
\label{fig:mesh-voro-agg}
\end{figure}

\begin{figure}
\centering
\caption{
Mesh obtained by aggregating horse cells. Boundaries of aggregates are marked in red.
}
\subfloat{\includegraphics[width=0.40\textwidth]{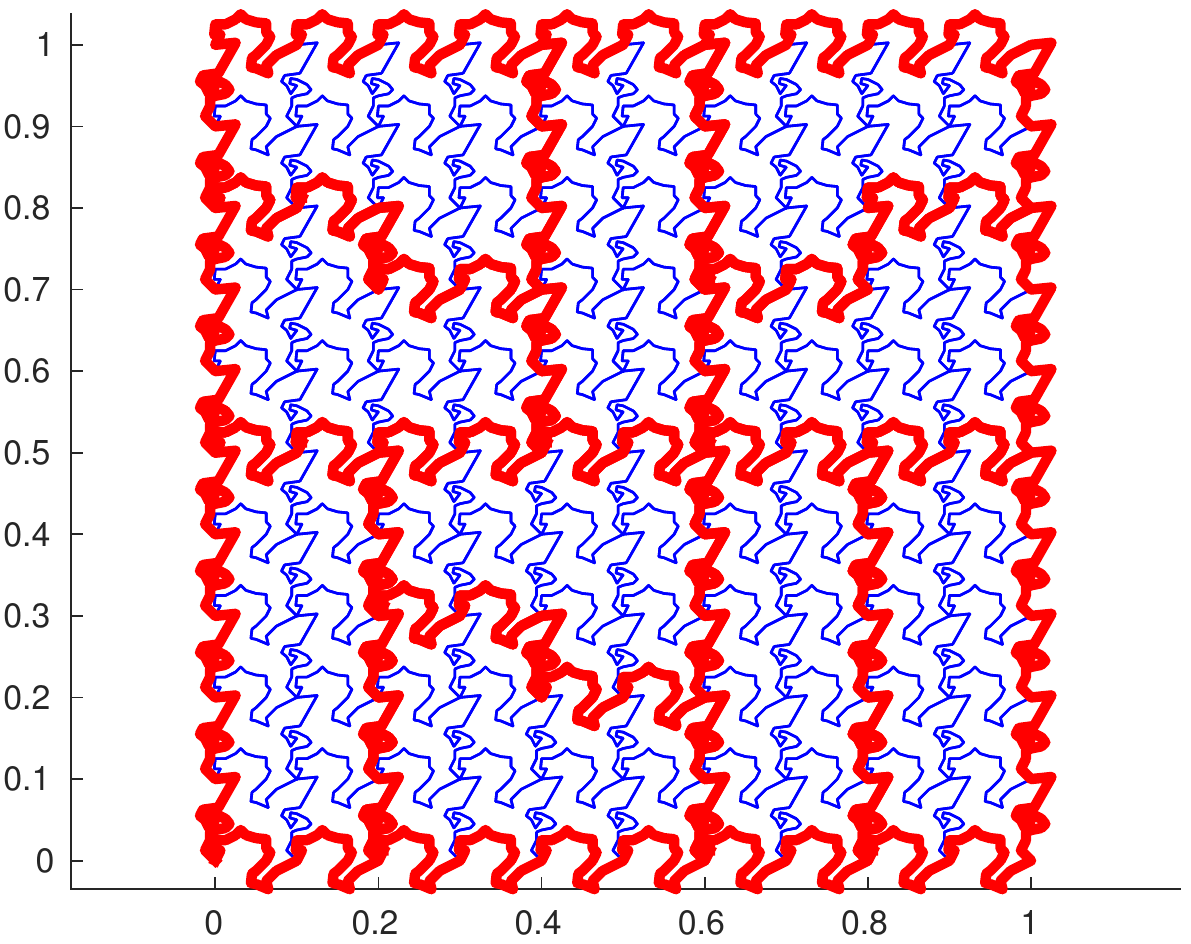}}
\label{fig:mesh-horse-agg}
\end{figure}

Moreover, in view of a possible use of AMG in an adaptive approach, we also deal with meshes obtained by agglomerating an underlying fine mesh. More precisely, 
we consider  different Voronoi and horse meshes and build coarse grid elements by agglomerating fine grid elements, 
that is, we deal with: 
\begin{enumerate}[1)]
\setcounter{enumi}{4}
	\item aggregates of Voronoi cells  (Figure~\ref{fig:mesh-voro-agg}), 
	\item 
aggregates of "horse" cells (Figure~\ref{fig:mesh-horse-agg}). 
\end{enumerate}
 
Let us now analyze the results.  The first Tables~\ref{tab:mesh-hexa},\ref{tab:mesh-voro},\ref{tab:mesh-horse},\ref{tab:mesh-voro-agg},\ref{tab:mesh-horse-agg} list the values of the geometrical parameters defined above for each mesh considered; 
 %
 we recall that the number of unknowns coincides with the number of vertices $N_\textup{v}$.
 Then,  in Tables~\ref{tab:k1-hexa},\ref{tab:k1-voro},\ref{tab:k1-horse},\ref{tab:k1-koch} and \ref{tab:k1-voro-agg},\ref{tab:k1-horse-agg} (for the case of aggregates of cells), 
we report, for each mesh, the condition number $\kappa$ 
of the matrix $A$ in (\ref{system}) with and without AMG preconditioning.
 The condition numbers are numerical approximations computed from the standard tridiagonal Lanczos matrix generated during the preconditioned CG iteration as the ratio between the maximum and the minimum eigenvalues, see e.g. \cite{GolubVanLoan_1996}.
Since in the following we deal also with a performance comparison of different (iterative and direct) solvers available in PETSc, 
the stopping criterion fixed for CG is defined by 
 the relative residual error 
  obtained with the direct solver SuperLU\_DIST (see \texttt{atol} reported in each Table).

We start by presenting the behavior of the condition number and of the iteration count when hexagonal meshes are taken into account.
In this case the scalability of all the AMG preconditioners considered is clearly shown in Table~\ref{tab:k1-hexa}: the iteration count to converge does not depend on the number of the mesh nodes. If we pass from hexagonal to  Voronoi meshes (see Table~\ref{tab:k1-voro}), which are characterized by worse geometrical parameters, we notice only a slight  increase in the number of iterations and in the condition number while the scalability is still preserved.

Instead, when dealing with horse or snowflake elements, we note  an increase in the number of iterations (Tables~\ref{tab:k1-horse}-\ref{tab:k1-koch}); scalability is not preserved by c-GAMG.
The best results are obtained by using BoomerAMG for the horse meshes and ML for the snowflake meshes. 
In this last case,  scalability is preserved also by BoomerAMG and a-GAMG, but with a higher number of iterations.
Similar results can also be observed with meshes of agglomerates of Voronoi or horse 
cells (see  Tables~\ref{tab:k1-voro-agg} and \ref{tab:k1-horse-agg}).  As before, BoomerAMG outperforms the other preconditioners, but, for the agglomerates of snowflakes cells, it still exhibits a high number of iterations, thereby worsening its computational performance, as we will see in the next sections. %




\begin{table}
	\begin{minipage}[b]{0.56\linewidth}
		\centering
			\caption{Meshes of regular hexagons; 
			for all these meshes, $\gamma_0 \approx \num[scientific-notation=fixed, fixed-exponent=0]{5.33}, \gamma_1 \approx \num[scientific-notation=fixed, fixed-exponent=0]{3.16}$.}
		\label{tab:mesh-hexa}
		\begin{tabular}{
				c
				S[table-format=7.0]
				S[table-format=7.0]
				S[table-format=1.{\roundPrecision}e-1]
				S[table-format=1.{\roundPrecision}e-1]
			}
			\toprule
			Mesh & {$N_\textup{elt}$} & {$N_\textup{v}$} & {$h$} & {$h_\textup{min}$}\\
			\midrule
			hexa$_1$ & 10151 & 20304 & 1.333333e-02 & 3.333333e-03\\
			hexa$_2$ & 40301 & 80604 & 6.666667e-03 & 1.666667e-03\\
			hexa$_3$ & 90451 & 180904 & 4.444444e-03 & 1.111111e-03\\
			hexa$_4$ & 160601 & 321204 & 3.333333e-03 & 8.333333e-04\\
			hexa$_5$ & 250751 & 501504 & 2.666667e-03 & 6.666667e-04\\
			hexa$_6$ & 360901 & 721804 & 2.222222e-03 & 5.555556e-04\\
			hexa$_7$ & 491051 & 982104 & 1.904762e-03 & 4.761905e-04\\
			hexa$_8$ & 641201 & 1282404 & 1.666667e-03 & 4.166667e-04\\
			hexa$_9$ & 811351 & 1622704 & 1.481481e-03 & 3.703704e-04\\
			hexa$_{10}$ & 1001501 & 2003004 & 1.333333e-03 & 3.333333e-04\\
			\bottomrule
		\end{tabular}
	\end{minipage}\hfill
\end{table}
\begin{table}
	\centering
	\caption{
	Condition number $\kappa$ of matrix $A$ 
		with and without AMG preconditioning on the {\em hexagonal meshes} of Table~\ref{tab:mesh-hexa}. 
	}
	\label{tab:k1-hexa}
	\small
	\begin{tabular}{c c c c c c c}
		\toprule
		\multirow{2}*{Mesh} & \multicolumn{1}{c}{\multirow{2}*{\texttt{rtol}}} & {No Prec} & {c-GAMG} & {a-GAMG} & {BoomerAMG} & {ML}\\
		\cmidrule(lr){3-7}
		&   &  {$\kappa$}   &   {$\kappa$}   &   {$\kappa$}   &   {$\kappa$}   &   {$\kappa$} (its) \\
		\midrule
		hexa$_{1}$   &   4.51e-16   &   7.74e+02  (109)  &   1.12  (10)  &   1.52  (15) &   1.31 (13) &   1.66 (15) \\
		hexa$_{2}$   &   9.36e-16   &   3.09e+03  (141)  &   1.14  (10)  &   1.75  (16) &   1.39 (14) &   1.97 (18) \\
		hexa$_{3}$   &   1.38e-15   &   6.96e+03  (194)  &   1.18  (11)  &   1.69  (17) &   1.49 (15) &   2.22 (20) \\
		hexa$_{4}$   &   1.85e-15   &   1.24e+04  (255)  &   1.18  (11)  &   1.92  (18) &   1.52 (15) &   2.51 (22) \\
		hexa$_{5}$   &   2.38e-15   &   1.93e+04  (316)  &   1.20  (11)  &   1.95  (18) &   1.56 (16) &   2.91 (23) \\
		hexa$_{6}$   &   2.80e-15   &   2.78e+04  (377)  &   1.23  (12)  &   2.01  (19) &   1.63 (16) &   3.12 (24) \\
		hexa$_{7}$   &   3.28e-15   &   3.79e+04  (438)  &   1.24  (12)  &   2.13  (19) &   1.61 (16) &   3.37 (25) \\
		hexa$_{8}$   &   3.73e-15   &   4.95e+04  (492)  &   1.26  (12)  &   2.18  (20) &   1.72 (17) &   3.87 (26) \\
		hexa$_{9}$   &   4.25e-15   &   6.26e+04  (545)  &   1.27  (12)  &   2.14  (20) &   1.72 (17) &   3.94 (26) \\
		hexa$_{10}$  &   4.89e-15   &   7.73e+04  (595)  &   1.27  (12)  &   2.25  (20) &   1.80 (17) &   4.05 (27) \\
		\bottomrule
	\end{tabular}
\end{table}



\begin{table}[h!]
\centering
\caption{Meshes of Voronoi cells used in the experiments.
}
\begin{tabular}
{
c
S[table-format=7.0]
S[table-format=7.0]
S[table-format=1.{\roundPrecision}e-1]
S[table-format=1.{\roundPrecision}e-1]
S[table-format=1.{\roundPrecision}e+1]
S[table-format=1.{\roundPrecision}e+1]
}
\toprule
Mesh & {$N_\textup{elt}$} & {$N_\textup{v}$} & {$h$} & {$h_\textup{min}$} & {$\gamma_0$} & {$\gamma_1$}\\
\midrule
voro$_1$ & 2500 & 5006 & 6.384666e-02 & 6.192528e-06 & 1.422215e+01 & 5.852731e+03\\
voro$_2$ & 5000 & 10008 & 4.344562e-02 & 5.845179e-07 & 1.456663e+01 & 3.404522e+04\\
voro$_3$ & 10000 & 20007 & 3.470002e-02 & 1.732139e-07 & 2.525382e+01 & 9.493668e+04\\
voro$_4$ & 20000 & 40011 & 2.405393e-02 & 2.138871e-07 & 2.087832e+01 & 7.246942e+04\\
voro$_5$ & 40000 & 80007 & 1.726980e-02 & 8.256465e-08 & 2.675024e+01 & 7.178744e+04\\
voro$_6$ & 80000 & 160028 & 1.140086e-02 & 5.998477e-09 & 2.881822e+01 & 1.100228e+06\\
voro$_7$ & 160000 & 320020 & 8.860795e-03 & 4.019971e-09 & 3.143276e+01 & 1.238004e+06\\
voro$_8$ & 320000 & 640035 & 6.248831e-03 & 3.592774e-09 & 3.074696e+01 & 9.566916e+05\\
voro$_9$ & 640000 & 1280053 & 4.361969e-03 & 1.139408e-09 & 4.041473e+01 & 1.914480e+06\\
voro$_{10}$ & 1280000 & 2560065 & 3.393600e-03 & 7.503513e-10 & 5.084927e+01 & 2.289585e+06\\
\bottomrule
\end{tabular}
\label{tab:mesh-voro}
\end{table}
\begin{table}[hb!]
	\centering
	\caption{
	Condition number $\kappa$ of matrix $A$ 
		with and without AMG preconditioning on the {\em Voronoi meshes} of Table~\ref{tab:mesh-voro}.
		With mesh voro$_{10}$, CG without preconditioning requires more than $10000$ iterations to reach convergence.}
	\label{tab:k1-voro}
	\small
	\begin{tabular}{c c c c c c c}
		\toprule
		\multirow{2}*{Mesh} & \multicolumn{1}{c}{\multirow{2}*{\texttt{rtol}}} & {No Prec} & {c-GAMG} & {a-GAMG} & {BoomerAMG} & {ML}\\
		\cmidrule(lr){3-7}
		&   &  {$\kappa$}   &   {$\kappa$}   &   {$\kappa$}   &   {$\kappa$}   &   {$\kappa$}\\
		\midrule
		voro$_{1}$   &   3.50e-14   &   2.79e+03  (371) &   1.43 (12)  &   2.09  (17) &  1.41  (13)  &   2.60 (19)\\
		voro$_{2}$   &   7.34e-14   &   6.91e+03  (553) &   1.63 (12)  &   2.50  (19) &  1.48  (13)  &   3.04 (22)\\
		voro$_{3}$   &   1.45e-13   &   1.55e+04  (827) &   1.81 (14)  &   2.65  (20) &  1.54  (13)  &   3.76 (24)\\
		voro$_{4}$   &   2.88e-13   &   3.44e+04 (1141) &   1.90 (15)  &   3.12  (21) &  1.55  (13)  &   4.29 (25)\\
		voro$_{5}$   &   5.83e-13   &   7.87e+04 (1736) &   2.04 (15)  &   3.27  (21) &  1.89  (14)  &   4.50 (27)\\
		voro$_{6}$   &   1.19e-12   &   1.36e+05 (2383) &   2.41 (17)  &   3.53  (22) &  1.80  (14)  &   5.56 (28)\\
		voro$_{7}$   &   2.41e-12   &   3.52e+05 (3633) &   2.72 (18)  &   3.99  (23) &  1.78  (14)  &   6.41 (30)\\
		voro$_{8}$   &   4.21e-12   &   8.98e+05 (5412) &   3.32 (19)  &   4.23  (23) &  2.03  (15)  &   8.20 (32)\\
		voro$_{9}$   &   9.70e-12   &   2.15e+06 (7775) &   3.92 (20)  &   4.77  (24) &  2.38  (15)  &   8.18 (33)\\
		voro$_{10}$  &   1.34e-11   &         {-}       &   4.11 (21)  &   5.24  (25) &  2.22  (15)  &  10.00 (35)\\
		\bottomrule
	\end{tabular}
\end{table}



\begin{table}
\centering
\caption{Meshes of horse cells  used in the experiments;
	for all of these meshes, $\gamma_1 = 39$.}
\begin{tabular}{
c
S[table-format=5.0]
S[table-format=6.0]
S[table-format=1.{\roundPrecision}e-1]
S[table-format=1.{\roundPrecision}e-1]
}
\toprule
Mesh & {$N_\textup{elt}$} & {$N_\textup{v}$} & {$h$} & {$h_\textup{min}$}\\
\midrule
horse$_1$ & 2500 & 94401 & 3.515412e-02 & 9.013878e-04\\
horse$_2$ & 3600 & 135481 & 2.929510e-02 & 7.511565e-04\\
horse$_3$ & 4900 & 183961 & 2.511009e-02 & 6.438484e-04\\
horse$_4$ & 6400 & 239841 & 2.197133e-02 & 5.633674e-04\\
horse$_5$ & 8100 & 303121 & 1.953007e-02 & 5.007710e-04\\
horse$_6$ & 10000 & 373801 & 1.757706e-02 & 4.506939e-04\\
horse$_7$ & 12100 & 451881 & 1.597915e-02 & 4.097217e-04\\
horse$_8$ & 14400 & 537361 & 1.464755e-02 & 3.755783e-04\\
horse$_9$ & 16900 & 630241 & 1.352082e-02 & 3.466876e-04\\
horse$_{10}$ & 19600 & 730521 & 1.255504e-02 & 3.219242e-04\\
\bottomrule
\end{tabular}
\label{tab:mesh-horse}
\end{table}

\begin{table}
	\centering
	\caption{
	Condition number $\kappa$ of matrix $A$ 
		with and without AMG preconditioning  on meshes of {\em horses} (see Table~\ref{tab:mesh-horse}).
	}
	\label{tab:k1-horse}
	\small
	\begin{tabular}{c c c c c c c}
		\toprule
		\multirow{2}*{Mesh} & \multicolumn{1}{c}{\multirow{2}*{\texttt{rtol}}} & {No Prec} & {c-GAMG} & {a-GAMG} & {BoomerAMG} & {ML}\\
		\cmidrule(lr){3-7}
		&   &  {$\kappa$}   &   {$\kappa$}   &   {$\kappa$}   &   {$\kappa$}   &   {$\kappa$}\\
		\midrule
		horse$_{1}$   &  1.15e-13 & 8.82e+04 (2071) & 127.55  (111) &  19.84 (59) &  3.42 (24) & 30.00 (72)\\
		horse$_{2}$   &  1.53e-13   &  1.28e+05  (2467) &   183.06  (130) &   20.63  (60) &   3.61  (24) &   30.43 (73)\\
		horse$_{3}$   &  1.95e-13   &  1.75e+05  (2837) &   243.36  (151) &   22.14  (62) &   3.49 (24) &   32.66 (76)\\
		horse$_{4}$   &  2.37e-13   &  2.29e+05  (3204) &   349.99  (176) &   22.87  (64) &   3.82  (25) &   36.46 (79)\\
		horse$_{5}$   &   2.83e-13  &  2.90e+05  (3583) &   535.99  (206) &   24.44  (66) &   3.65  (25) &   38.29 (81)\\
		horse$_{6}$   &   3.35e-13  &  3.60e+05  (3953) &   485.97  (212) &   25.30  (67) &   3.81 (25) &   41.80 (83)\\
		horse$_{7}$   &   3.85e-13  &  4.36e+05  (4273) &   641.35  (232) &   25.57  (67) &   4.03  (25) &   41.53 (84)\\
		horse$_{8}$   &   4.41e-13  &  5.19e+05  (4644) &   741.34  (261) &   26.23  (68) &   4.05  (25) &   39.89 (85)\\
		horse$_{9}$   &   4.98e-13  &  6.10e+05  (5007) &   894.57  (280) &   27.14 (69)  &   4.00  (25) &   41.65 (84)\\
		horse$_{10}$  &   5.57e-13  &  7.08e+05  (5375) &   1009.26 (298) &   27.59  (69) &   4.11  (25) &   46.51 (89)\\
		\bottomrule
	\end{tabular}
\end{table}

%

\begin{table}
\centering
\caption{Meshes with snowflakes cells. 
For all of these meshes, $\gamma_1 \approx \num[scientific-notation=fixed, fixed-exponent=0]{6.854545e+01}$.}
\begin{tabular}{
c
S[table-format=4.0]
S[table-format=6.0]
S[table-format=1.{\roundPrecision}e-1]
S[table-format=1.{\roundPrecision}e-1]
}
\toprule
Mesh & {$N_\textup{elt}$} & {$N_\textup{v}$} & {$h$} & {$h_\textup{min}$}\\
\midrule
koch$_1$ & 288 & 27973 & 1.178511e-01 & 1.719314e-03\\
koch$_2$ & 512 & 49713 & 8.838835e-02 & 1.289485e-03\\
koch$_3$ & 800 & 77661 & 7.071068e-02 & 1.031588e-03 \\
koch$_4$ & 1152 & 111817 & 5.892557e-02 & 8.596569e-04\\
koch$_5$ & 1568 & 152181 & 5.050763e-02 & 7.368488e-04\\
koch$_6$ & 2048 & 198753 & 4.419417e-02 & 6.447427e-04\\
koch$_7$ & 2592 & 251533 & 3.928371e-02 & 5.731046e-04\\
koch$_8$ & 3200 & 310521 & 3.535534e-02 & 5.157941e-04\\
koch$_9$ & 3872 & 375717 & 3.214122e-02 & 4.689038e-04\\
koch$_{10}$ & 4608 & 447121 & 2.946278e-02 & 4.298285e-04\\
koch$_{11}$ & 5408 & 524733 & 2.719641e-02 & 3.967647e-04\\
koch$_{12}$ & 6272 & 608553 & 2.525381e-02 & 3.684244e-04\\
koch$_{13}$ & 7200 & 698581 & 2.357023e-02 & 3.438628e-04\\
\bottomrule
\end{tabular}
\label{tab:mesh-koch}
\end{table}
\begin{table}[htb]
	\centering
	\caption{
	Condition number $\kappa$ of matrix $A$ 
		on meshes with {\em  snowflakes}. 
	}
	\label{tab:k1-koch}
	\small
	\begin{tabular}{c c c c c c c}
		\toprule
		\multirow{2}*{Mesh} & \multicolumn{1}{c}{\multirow{2}*{\texttt{rtol}}} & {No Prec} & {c-GAMG} & {a-GAMG} & {BoomerAMG} & {ML}\\
		\cmidrule(lr){3-7}
		&   &  {$\kappa$}   &   {$\kappa$}   &   {$\kappa$}   &   {$\kappa$}   &   {$\kappa$}\\
		\midrule
		koch$_{1}$    &   3.87e-13   &   8.41e+04 (1069) &   46.51  (50)  &   15.49  (50) &   64.08  (91) &   8.76 (35)\\
		koch$_{2}$    &   7.38e-13  &   1.51e+05   (1430)  &   79.63  (63)  &   16.47 (53)  &   68.38   (91)  &   9.15 (35)\\
		koch$_{3}$    &   1.21e-12    &    2.37e+05  (1772)  &   126.33 (77)  &   17.38 (54)  &   71.51  (94)  &   10.40 (35)\\
		koch$_{4}$    &   1.74e-12   &   3.42e+05   (2117) &   181.52  (89)  &   18.25  (55) &   73.83  (94)  &   12.98 (35)\\
		koch$_{5}$    &   2.43e-12   &   4.67e+05  (2464)  &   251.53  (102) &   20.44 (56)  &   75.60   (93)  &   13.80 (35)\\
		koch$_{6}$    &   3.19e-12  &   6.10e+05   (2632) &   302.75  (110)  &   19.16 (55)    &   76.98  (92)  &   15.73 (35)\\
		koch$_{7}$    &   4.01e-12  &   7.73e+05  (2930)  &   414.78  (127)  &   19.64  (55)  &   77.81  (93) &   15.87 (35)\\
		koch$_{8}$    &   4.96e-12  &   9.55e+05  (3182)  &   507.17  (138)  &   20.36 (56)  &   79.03   (93) &   16.05 (35)\\
		koch$_{9}$    &   6.07e-12  &   1.16e+06   (3469)   &   629.20  (152)  &   20.06 (54)  &   79.80   (92)  &   17.63 (35)\\
		koch$_{10}$   &   7.38e-12  &   1.38e+06  (3670)  &   724.47  (158)  &   19.32 (54)  &   77.60    (91) &   18.58 (35)\\
		koch$_{11}$   &   8.66e-12  &   1.61e+06  (3915)   &   867.05   (170) &   19.76  (54)  &   81.06  (91)  &   19.66 (35)\\
		koch$_{12}$   &   9.74e-12  &   1.87e+06 (4141)   &   1018.24  (182)  &   20.40  (54)  &   81.51  (91)  &   19.89 (35)\\
		koch$_{13}$   &   1.14e-11    &   2.15e+06  (4384)   &   1101.59  (187)  &   19.97 (54)   &   80.98  (91)  &   20.52 (35)\\
		\bottomrule
	\end{tabular}
\end{table}



\begin{table}
\centering
\caption{Meshes of aggregates of Voronoi cells used in the experiments.
 }
\begin{tabular}{
c
S[table-format=6.0]
S[table-format=7.0]
S[table-format=1.{\roundPrecision}e-1]
S[table-format=1.{\roundPrecision}e-1]
S[table-format=1.{\roundPrecision}e+1]
}
\toprule
Mesh & {$N_\textup{elt}$} & {$N_\textup{v}$} & {$h$} & {$h_\textup{min}$} & {$\gamma_1$}\\
\midrule
a-voro$_1$ & 512 & 5652 & 1.261054e-01 & 7.974791e-07 & 9.732997e+04\\
a-voro$_2$ & 1024 & 11190 & 9.861447e-02 & 1.732139e-07 & 3.051253e+05\\
a-voro$_3$ & 2048 & 22250 & 6.306131e-02 & 3.470813e-07 & 1.148450e+05\\
a-voro$_4$ & 4096 & 44615 & 4.800332e-02 & 8.256465e-08 & 3.307747e+05\\
a-voro$_5$ & 8192 & 89283 & 3.770308e-02 & 2.438264e-08 & 5.740921e+05\\
a-voro$_6$ & 16384 & 180219 & 2.780324e-02 & 4.019971e-09 & 4.386855e+06\\
a-voro$_7$ & 32768 & 365921 & 1.879699e-02 & 3.592774e-09 & 2.347532e+06\\
a-voro$_8$ & 65536 & 731552 & 1.399502e-02 & 1.139408e-09 & 6.975626e+06\\
a-voro$_9$ & 131072 & 1471656 & 1.023166e-02 & 3.431785e-09 & 1.411328e+06\\
\bottomrule
\end{tabular}
\label{tab:mesh-voro-agg}
\end{table}
\begin{table}
	\centering
	\caption{{\em Aggregates of Voronoi cells} (Table~\ref{tab:mesh-voro-agg}). Condition number $\kappa$ of matrix $A$ 
		with and without AMG preconditioning. 
	}
	\label{tab:k1-voro-agg}
	\small
	\begin{tabular}{c c c c c c c}
		\toprule
		\multirow{2}*{Mesh} & \multicolumn{1}{c}{\multirow{2}*{\texttt{rtol}}} & {No Prec} & {c-GAMG} & {a-GAMG} & {BoomerAMG} & {ML}\\
		\cmidrule(lr){3-7}
		&   &  {$\kappa$}   &   {$\kappa$}   &   {$\kappa$}   &   {$\kappa$}   &   {$\kappa$}\\
		\midrule
		a-voro$_{1}$   &   4.20e-14   & 2.89e+03  (380)  &   3.45 (19)  &  4.00 (24) &  1.92 (16) &   5.37 (29)\\
		a-voro$_{2}$   &   8.10e-14   & 1.06e+04  (586)  &   5.68 (25)  &  6.01 (30) &  2.09 (17) &   8.34 (35)\\
		a-voro$_{3}$   &   1.60e-13   & 1.62e+04  (833)  &   8.87 (31)  &  6.73 (32) &  1.90 (16) &   8.62 (37)\\
		a-voro$_{4}$   &   3.24e-13   & 3.84e+04 (1258)  &  14.49 (41)  &  7.67 (35) &  2.06 (16) &  10.33 (41)\\
		a-voro$_{5}$   &   6.67e-13   & 1.04e+05 (1817)  &  28.81 (55)  & 10.28 (40) &  2.51 (17) &  14.89 (47)\\
		a-voro$_{6}$   &   1.37e-12   & 2.03e+05 (2710)  &  79.60 (77)  & 10.96 (41) &  4.14 (19) &  17.58 (49)\\
		a-voro$_{7}$   &   2.81e-12   & 4.12e+05 (3994)  & 141.87 (114) & 12.42 (42) &  2.62 (18) &  19.46 (52)\\
		a-voro$_{8}$   &   5.70e-12   & 9.36e+05 (5690)  & 247.73 (154) & 13.56 (44) &  2.85 (18) &  23.09 (56)\\
		a-voro$_{9}$   &   8.55e-13   & 2.05e+06 (8410)  & 515.85 (229) & 14.94 (49) &  3.22 (20) &  24.07 (63)\\
		\bottomrule
	\end{tabular}
\end{table}



\begin{table}
\centering
\caption{Meshes of aggregated horse cells used in the experiments.
}
\begin{tabular}{
c
S[table-format=4.0]
S[table-format=6.0]
S[table-format=1.{\roundPrecision}e-1]
S[table-format=1.{\roundPrecision}e-1]
S[table-format=1.{\roundPrecision}e+1]
}
\toprule
Mesh & {$N_\textup{elt}$} & {$N_\textup{v}$} & {$h$} & {$h_\textup{min}$} & {$\gamma_1$}\\
\midrule
a-horse$_1$ & 250 & 42621 & 1.614181e-01 & 9.013878e-04 & 1.790773e+02\\
a-horse$_2$ & 360 & 61151 & 1.471153e-01 & 7.511565e-04 & 1.958518e+02\\
a-horse$_3$ & 490 & 82007 & 1.250017e-01 & 6.438484e-04 & 1.941476e+02\\
a-horse$_4$ & 640 & 107899 & 1.080152e-01 & 5.633674e-04 & 1.917314e+02\\
a-horse$_5$ & 810 & 134589 & 1.003632e-01 & 5.007710e-04 & 2.004174e+02\\
a-horse$_6$ & 1000 & 168525 & 9.032692e-02 & 4.506939e-04 & 2.004174e+02\\
\bottomrule
\end{tabular}
\label{tab:mesh-horse-agg}
\end{table}
\begin{table}
	\centering
	\caption{
		 Condition number $\kappa$ of matrix $A$ 
		with and without AMG preconditioning on meshes of {\em aggregated horses}.
		ML fails on all the meshes due to an internal \emph{Segment Violation} error.}
	\label{tab:k1-horse-agg}
	\small
	\begin{tabular}{c c c c c c c}
		\toprule
		\multirow{2}*{Mesh} & \multicolumn{1}{c}{\multirow{2}*{\texttt{rtol}}} & {No Prec} & {c-GAMG} & {a-GAMG} & {BoomerAMG} & {ML}\\
		\cmidrule(lr){3-7}
		&   &  {$\kappa$}   &   {$\kappa$}   &   {$\kappa$}   &   {$\kappa$}   &   {$\kappa$}\\
		\midrule
		a-horse$_{1}$   &   8.31e-14   &   1.44e+05 (1896) &  182.00 (82)   &   51.36 (73)  &  22.66  (60)  &   {-}\\
		a-horse$_{2}$   &   1.28e-13   &   2.58e+05 (2625) &  267.76 (106)  &   61.17 (87)  &  40.17  (69)  &   {-}\\
		a-horse$_{3}$   &   1.81e-13   &   3.37e+05 (3221) &  355.49 (125)  &   67.09 (97)  &  41.09  (69)  &   {-}\\
		a-horse$_{4}$   &   1.95e-13   &   4.61e+05 (3874) &  473.41 (151)  &   78.85 (108) &  25.67  (69)  &   {-}\\
		a-horse$_{5}$   &   2.71e-13   &   6.47e+05 (4456) &  643.31 (171)  &   82.15 (115) &  34.73  (67)  &   {-}\\
		a-horse$_{6}$   &   2.83e-13   &   8.11e+05 (5099) &  728.71 (189)  &   96.54 (120) &  39.74  (75)  &   {-}\\
		\bottomrule
	\end{tabular}
\end{table}



 \subsection{Discontinuous coefficients}
 \label{sec:lowest-random}
We now want to test robustness of AMG preconditioners when discontinuous, highly heterogeneous coefficients $\rho$ are considered. 
We build a ``checkerboard'' pattern as follows:
we deal with five of the hexagonal meshes in Table~\ref{tab:mesh-hexa} (in particular, hexa$_2$, hexa$_4$, hexa$_6$, hexa$_8$ and hexa$_{10}$) and the eight 
 finest meshes in Table~\ref{tab:k1-voro} and partition their polygons into an increasing number of parts $L$ ($L = 64, 128, 256, 512$ and $1024$) by using METIS~\cite{metis-web-page}. 
Every element in a given part is then assigned the same diffusion coefficient $\rho = 10^\alpha$, with $\alpha$ random integer in $[-5, 5]$. 
 The loading term is random uniform in $[-1, 1]$.

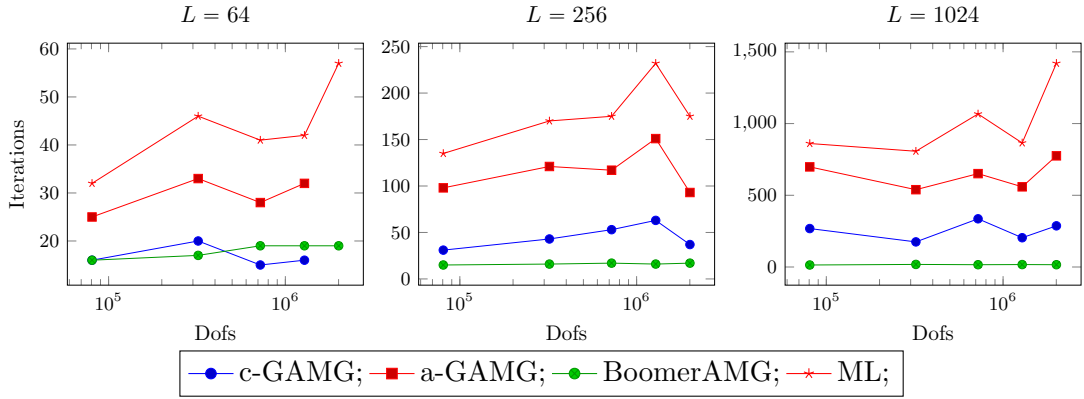
\begin{figure}
	\centering
	\begin{tikzpicture}[scale=0.8][trim axis right]
	\begin{semilogxaxis}[small, xlabel = {Dofs}, title = {$L = 64$},  ylabel={Iterations}, 
	legend columns=4,
	legend entries={$\text{c-GAMG}$;, $\text{a-GAMG}$;, $\text{BoomerAMG}$;, $\text{ML}$;},
	legend to name=robust-irho3-hexa,
	cycle list name=mycolorlist]
	\addplot coordinates {
		(80604, 16)
		(321204, 20)
		(721804, 15)
		(1282404, 16)
	};
	\addplot coordinates {
		(80604, 25)
		(321204, 33)
		(721804, 28)
		(1282404, 32)
	};
	\addplot coordinates {
		(80604, 16)
		(321204, 17)
		(721804, 19)
		(1282404, 19)
		(2003004, 19)
	};
	\addplot coordinates {
		(80604, 32)
		(321204, 46)
		(721804, 41)
		(1282404, 42)
		(2003004, 57)
	};
	\end{semilogxaxis}
	\end{tikzpicture}
	\begin{tikzpicture}[scale=0.8][trim axis right]
	\begin{semilogxaxis}[small, xlabel = {Dofs}, title = {$L = 256$}, 
	cycle list name=mycolorlist]
	\addplot coordinates {
		(80604, 31)
		(321204, 43)
		(721804, 53)
		(1282404, 63)
		(2003004, 37)
	};
	\addplot coordinates {
		(80604, 98)
		(321204, 121)
		(721804, 117)
		(1282404, 151)
		(2003004, 93)
	};
	\addplot coordinates {
		(80604, 15)
		(321204, 16)
		(721804, 17)
		(1282404, 16)
		(2003004, 17)
	};
	\addplot coordinates {
		(80604, 135)
		(321204, 170)
		(721804, 175)
		(1282404, 232)
		(2003004, 175)
	};
	\end{semilogxaxis}
	\end{tikzpicture}
	\begin{tikzpicture}[scale=0.8][trim axis right]
	\begin{semilogxaxis}[small, xlabel = {Dofs}, title = {$L = 1024$},  
	cycle list name=mycolorlist]
	\addplot coordinates {
		(80604, 268)
		(321204, 175)
		(721804, 336)
		(1282404, 204)
		(2003004, 287)
	};
	\addplot coordinates {
		(80604, 698)
		(321204, 539)
		(721804, 651)
		(1282404, 559)
		(2003004, 775)
	};
	\addplot coordinates {
		(80604, 14)
		(321204, 18)
		(721804, 16)
		(1282404, 17)
		(2003004, 16)
	};
	\addplot coordinates {
		(80604, 861)
		(321204, 807)
		(721804, 1067)
		(1282404, 865)
		(2003004, 1420)
	};
	\end{semilogxaxis}
	\end{tikzpicture}
	\ref{robust-irho3-hexa}
	\caption{Discontinuous coefficients $\rho$. Number of iterations of PCG  with different AMG preconditioners in a  logarithmic scale for the x-axis and for the hexagonal meshes hexa$_2$, hexa$_4$, hexa$_6$, hexa$_8$ and hexa$_{10}$ of  Table~\ref{tab:mesh-hexa} partitioned 
	in $L = 64, 256, 1024$ parts.}
	\label{fig:hexa-random-amg}
\end{figure}

\begin{figure}
	\centering
	\begin{tikzpicture}[scale=0.8][trim axis right]
	\begin{semilogxaxis}[small,   xlabel=Dofs, title={ $L = 64$},  ylabel={Iterations}, 
	legend columns=4,
	legend entries={$\text{c-GAMG}$;, $\text{a-GAMG}$;, $\text{BoomerAMG}$;, $\text{ML}$;},
	legend to name=robust-irho3-voro,
	cycle list name=mycolorlist]
	%
	\addplot coordinates {
		(20007, 26)
		(40011, 30)
		(80007, 27)
		(160028, 29)
		(320020, 29)
		(640035, 38)
		(1280053, 38)
		(2560065, 50)
	};
	\addplot coordinates {
		(20007, 38)
		(40011, 38)
		(80007, 34)
		(160028, 40)
		(320020, 37)
		(640035, 57)
		(1280053, 47)
		(2560065, 52)
	};
	\addplot coordinates {
		(20007, 19)
		(40011, 20)
		(80007, 18)
		(160028, 19)
		(320020, 21)
		(640035, 20)
		(1280053, 22)
		(2560065, 23)
	};
	\addplot coordinates {
		(20007, 51)
		(40011, 55)
		(80007, 54)
		(160028, 49)
		(320020, 55)
		(640035, 70)
		(1280053, 71)
		(2560065, 89)
	};
	\end{semilogxaxis}
	\end{tikzpicture}
	\begin{tikzpicture}[scale=0.8][trim axis right]
	\begin{semilogxaxis}[small,   xlabel=Dofs, title={ $L = 256$}, 
	cycle list name=mycolorlist]
	\addplot coordinates {
		(20007, 110)
		(40011, 79)
		(80007, 135)
		(160028, 150)
		(320020, 178)
		(640035, 151)
		(1280053, 195)
		(2560065, 250)
	};
	\addplot coordinates {
		(20007, 198)
		(40011, 115)
		(80007, 163)
		(160028, 193)
		(320020, 187)
		(640035, 152)
		(1280053, 248)
		(2560065, 282)
	};
	\addplot coordinates {
		(20007, 15)
		(40011, 18)
		(80007, 17)
		(160028, 17)
		(320020, 18)
		(640035, 19)
		(1280053, 19)
		(2560065, 21)
	};
	\addplot coordinates {
		(20007, 249)
		(40011, 158)
		(80007, 268)
		(160028, 313)
		(320020, 288)
		(640035, 251)
		(1280053, 381)
		(2560065, 387)
	};
	\end{semilogxaxis}
	\end{tikzpicture}
	\begin{tikzpicture}[scale=0.8][trim axis right]
	\begin{semilogxaxis}[small, xlabel=Dofs,  title={ $L = 1024$}, 
	 cycle list name=mycolorlist]
	\addplot coordinates {
		(20007, 693)
		(40011, 836)
		(80007, 683)
		(160028, 984)
		(320020, 735)
		(640035, 954)
		(1280053, 1018)
		(2560065, 878)
	};
	\addplot coordinates {
		(20007, 1078)
		(40011, 1155)
		(80007, 1056)
		(160028, 1192)
		(320020, 960)
		(640035, 1076)
		(1280053, 1149)
		(2560065, 867)
	};
	\addplot coordinates {
		(20007, 16)
		(40011, 17)
		(80007, 17)
		(160028, 17)
		(320020, 18)
		(640035, 19)
		(1280053, 19)
		(2560065, 21)
	};
	\addplot coordinates {
		(20007, 1439)
		(40011, 1687)
		(80007, 1491)
		(160028, 1756)
		(320020, 1618)
		(640035, 2027)
		(1280053, 2078)
		(2560065, 1554)
	};
	\end{semilogxaxis}
	\end{tikzpicture}
	\ref{robust-irho3-voro}
	\caption{Discontinuous coefficients $\rho$. Number of iterations of PCG  with different AMG preconditioners 
		in a logarithmic scale for the x-axis and for the Voronoi meshes voro$_3$--voro$_{10}$ of Table~\ref{tab:k1-voro} partitioned 
		in $L = 64, 256, 1024$ parts. }
	%
	\label{fig:voro-random-metis-amg}
\end{figure}
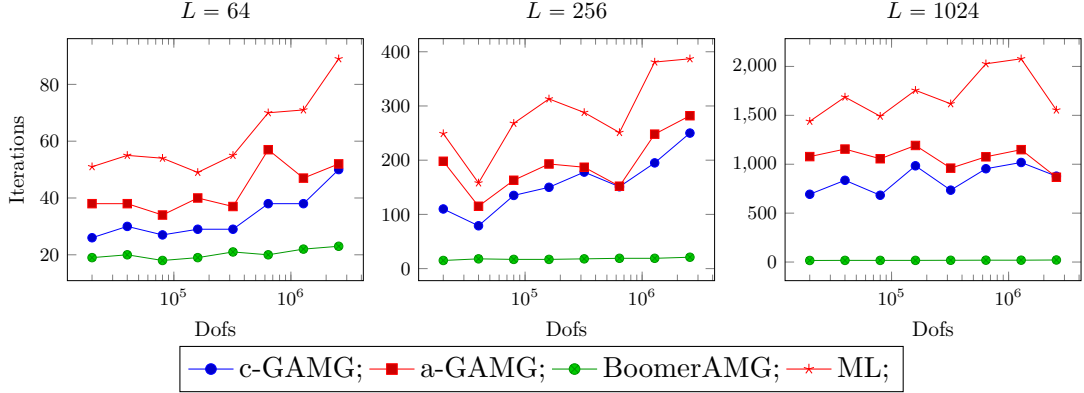

Here, convergence is decided by the \emph{absolute} size of the residual norm \texttt{abstol} obtained with the direct solver SuperLU\_DIST (see Tables~\ref{tab:abstol-hexa},~\ref{tab:abstol-voro}). Results are displayed  in Table~\ref{tab:hexa-random} and Figure~\ref{fig:hexa-random-amg} for hexagonal meshes, and in Table~\ref{tab:voro-random-metis} and Figure~\ref{fig:voro-random-metis-amg} for Voronoi meshes. 
We point out that, with highly oscillating coefficients, CG without AMG preconditioning does not attain convergence on any mesh.

Tables~\ref{tab:hexa-random} and \ref{tab:voro-random-metis} show that only BoomerAMG, for both hexagonal and Voronoi meshes, is performing well even with highly oscillating coefficients; the iteration count of BoomerAMG  does not 
depend  on the size of the jump of the diffusion coefficient and on the problem size,  see Figures ~\ref{fig:hexa-random-amg}-\ref{fig:voro-random-metis-amg} where we used a logarithmic scale for the x-axis.

\begin{table}
\centering
\caption{Absolute size of the residual norm \texttt{abstol} obtained with the direct solver SuperLU\_DIST on the hexagonal meshes partitioned with METIS with random data.}
\label{tab:abstol-hexa}
\begin{tabular}
{
c
S[table-format=1.{\roundPrecision}e-2]
S[table-format=1.{\roundPrecision}e-2]
S[table-format=1.{\roundPrecision}e-2]
S[table-format=1.{\roundPrecision}e-2]
S[table-format=1.{\roundPrecision}e-2]
}
\toprule
$\text{Mesh} \backslash L$ & {$64$} & {$128$} & {$256$} & {$512$} & {$1024$}\\
\midrule
hexa$_2$   &   3.377076e-15   &   2.749819e-12   &   1.659726e-13   &   6.289179e-13   &   1.097806e-12\\
hexa$_4$   &   4.656541e-15   &   3.928083e-14   &   3.005561e-13   &   1.402578e-13   &   1.724394e-13\\
hexa$_6$   &   7.400275e-16   &   6.396524e-15   &   3.228314e-13   &   9.881418e-13   &   1.651562e-12\\
hexa$_8$   &   1.146216e-15   &   7.971567e-13   &   4.437052e-13   &   4.136265e-13   &   1.156183e-13\\
hexa$_{10}$   &   2.303510e-15   &   1.376253e-12   &   1.467877e-13   &   1.327990e-13   &   4.095764e-12\\
\bottomrule
\end{tabular}
\end{table}

\begin{table}
\centering
\caption{Absolute size of the residual norm \texttt{abstol} obtained with the direct solver SuperLU\_DIST on the Voronoi meshes partitioned with METIS with random data.}
\label{tab:abstol-voro}
\begin{tabular}
{
c
S[table-format=1.{\roundPrecision}e-2]
S[table-format=1.{\roundPrecision}e-2]
S[table-format=1.{\roundPrecision}e-2]
S[table-format=1.{\roundPrecision}e-2]
S[table-format=1.{\roundPrecision}e-2]
}
\toprule
$\text{Mesh} \backslash L$ & {$64$} & {$128$} & {$256$} & {$512$} & {$1024$}\\
\midrule
voro$_{3}$   &   3.643422e-16   &   3.265309e-12   &   4.943836e-13   &   1.260445e-12   &   2.146567e-13\\
voro$_{4}$   &   4.640741e-15   &   4.042216e-14   &   2.189554e-14   &   8.516365e-14   &   3.836524e-13\\
voro$_{5}$   &   1.123398e-14   &   6.368778e-13   &   4.922339e-13   &   5.050658e-13   &   2.665815e-13\\
voro$_{6}$   &   5.361608e-15   &   1.154381e-12   &   1.684534e-13   &   3.386050e-12   &   5.315959e-13\\
voro$_{7}$   &   7.146989e-16   &   4.619196e-13   &   3.049112e-13   &   1.468307e-13   &   1.439853e-12\\
voro$_{8}$   &   1.287743e-14   &   3.819394e-14   &   5.987926e-14   &   1.578116e-13   &   1.529344e-12\\
voro$_{9}$   &   9.557891e-16   &   1.146418e-12   &   3.469132e-13   &   1.830300e-12   &   3.798424e-13\\
voro$_{10}$   &   2.582220e-15   &   3.416954e-14   &   3.007990e-13   &   3.142518e-13   &   1.439913e-13\\
\bottomrule
\end{tabular}
\end{table}

\begin{table}
\centering
\caption{Discontinuous coefficients $\rho$. Condition number $\kappa$ of matrix $A$ 
	with AMG preconditioning on the hexagonal meshes hexa$_2$, hexa$_4$, hexa$_6$, hexa$_8$ and hexa$_{10}$ of Table~\ref{tab:mesh-hexa} partitioned in $L$ parts. }
\label{tab:hexa-random}
\small
\begin{tabular}{
c
S[table-format=1.{\roundPrecision}e+1]
S[table-format=1.{\roundPrecision}e+1]
S[table-format=1.{\roundPrecision}e+1]
S[table-format=1.{\roundPrecision}e+1]
S[table-format=1.{\roundPrecision}e+1]
}
\toprule
\multicolumn{6}{c}{{c-GAMG}}\\
\midrule
$\text{Mesh} \backslash L$ & {$64$} & {$128$} & {$256$} & {$512$} & {$1024$}\\
\midrule
hexa$_{2}$   &   2.297887e+01   &   2.425112e+03   &   1.252988e+02   &   2.248741e+03   &   6.383563e+03\\
hexa$_{4}$   &   6.477467e+01   &   2.105721e+02   &   4.076302e+03   &   1.601431e+03   &   1.721920e+04\\
hexa$_{6}$   &   1.851394e+00   &   6.693057e+01   &   1.871567e+04   &   6.042888e+03   &   4.731036e+04\\
hexa$_{8}$   &   4.729968e+00   &   {-}   &   4.948828e+03   &   3.070528e+04   &   2.101103e+04\\
hexa$_{10}$  &   {-}   &   1.358146e+03   &   9.636772e+02   &   2.215524e+03   &   8.505207e+04\\
\midrule
\multicolumn{6}{c}{{a-GAMG}}\\
\midrule
$\text{Mesh} \backslash L$ & {$64$} & {$128$} & {$256$} & {$512$} & {$1024$}\\
\midrule
hexa$_{2}$   &   2.159164e+01   &   3.254016e+05   &   1.212575e+04   &   2.753859e+05   &   1.227976e+06\\
hexa$_{4}$   &   6.943918e+02   &   8.058224e+02   &   1.114240e+04   &   6.566179e+03   &   5.091468e+04\\
hexa$_{6}$   &   2.212608e+01   &   2.281328e+02   &   1.003750e+04   &   4.958900e+04   &   5.056637e+05\\
hexa$_{8}$   &   1.593507e+02   &   {-}   &   4.055078e+04   &   2.071823e+05   &   4.212288e+04\\
hexa$_{10}$  &   {-}   &   1.010984e+05   &   2.985231e+03   &   2.985851e+04   &   3.955145e+05\\
\midrule
\multicolumn{6}{c}{{BoomerAMG}}\\
\midrule
$\text{Mesh} \backslash L$ & {$64$} & {$128$} & {$256$} & {$512$} & {$1024$}\\
\midrule
hexa$_{2}$   &   1.677431e+00   &   1.969824e+00   &   1.778734e+00   &   2.122178e+00   &   1.848879e+00\\
hexa$_{4}$   &   1.899265e+00   &   1.787248e+00   &   1.797164e+00   &   1.875232e+00   &   2.555856e+00\\
hexa$_{6}$   &   2.688183e+00   &   2.023706e+00   &   2.035363e+00   &   2.521954e+00   &   2.140209e+00\\
hexa$_{8}$   &   2.007103e+00   &   2.033311e+00   &   1.955653e+00   &   2.364937e+00   &   2.076121e+00\\
hexa$_{10}$   &   1.901463e+00   &   2.066522e+00   &   2.058068e+00   &   2.158301e+00   &   2.658213e+00\\
\midrule
\multicolumn{6}{c}{{ML}}\\
\midrule
$\text{Mesh} \backslash L$ & {$64$} & {$128$} & {$256$} & {$512$} & {$1024$}\\
\midrule
hexa$_{2}$   &   2.850350e+02   &   1.315568e+06   &   1.220754e+04   &   4.016915e+05   &   7.808395e+05\\
hexa$_{4}$   &   1.558677e+03   &   3.600116e+03   &   2.629802e+04   &   1.704398e+05   &   7.012014e+04\\
hexa$_{6}$   &   1.135980e+02   &   1.028641e+03   &   9.517323e+03   &   1.666517e+05   &   1.206052e+06\\
hexa$_{8}$   &   1.073052e+03   &   1.127494e+05   &   1.121612e+05   &   1.692994e+06   &   6.880973e+05\\
hexa$_{10}$   &   2.114068e+02   &   4.810950e+05   &   1.017082e+04   &   7.783354e+05   &   6.844360e+05\\
\bottomrule
\end{tabular}
\end{table}

\begin{table}
\centering
\caption{Discontinuous coefficients $\rho$. Condition number $\kappa$ of matrix $A$ 
	with AMG preconditioning on the Voronoi meshes voro$_3$--voro$_{10}$ of 
	Table~\ref{tab:k1-voro} partitioned in $L$ parts. } 
\label{tab:voro-random-metis}
\small
\begin{tabular}{
c
S[table-format=1.{\roundPrecision}e+1]
S[table-format=1.{\roundPrecision}e+1]
S[table-format=1.{\roundPrecision}e+1]
S[table-format=1.{\roundPrecision}e+1]
S[table-format=1.{\roundPrecision}e+1]
}
\toprule
\multicolumn{6}{c}{{c-GAMG}}\\
\midrule
$\text{Mesh} \backslash L$ & {$64$} & {$128$} & {$256$} & {$512$} & {$1024$}\\
\midrule
voro$_{3}$   &   1.954338e+01   &   4.810408e+05   &   8.472994e+03   &   1.625234e+05   &   6.120649e+05\\
voro$_{4}$   &   2.988849e+02   &   7.448958e+02   &   1.967296e+03   &   2.162731e+04   &   5.143102e+04\\
voro$_{5}$   &   4.544447e+02   &   8.709399e+03   &   1.725037e+04   &   2.949029e+05   &   1.169738e+05\\
voro$_{6}$   &   6.249303e+01   &   2.290299e+04   &   2.088721e+04   &   6.278002e+03   &   3.067088e+05\\
voro$_{7}$   &   3.663888e+01   &   1.021271e+06   &   2.118504e+04   &   2.431300e+05   &   8.012810e+05\\
voro$_{8}$   &   1.096859e+02   &   7.894750e+03   &   2.999480e+04   &   8.593241e+04   &   4.306101e+05\\
voro$_{9}$   &   7.711819e+01   &   2.838958e+05   &   1.864640e+04   &   3.187226e+05   &   2.595585e+05\\
voro$_{10}$   &   8.500128e+02   &   1.402728e+04   &   1.138845e+05   &   9.281615e+04   &   6.332722e+05\\
\midrule
\multicolumn{6}{c}{{a-GAMG}}\\
\midrule
$\text{Mesh} \backslash L$ & {$64$} & {$128$} & {$256$} & {$512$} & {$1024$}\\
\midrule
voro$_{3}$   &   2.147681e+01   &   8.809488e+05   &   3.447852e+04   &   1.777274e+05   &   5.598563e+05\\
voro$_{4}$   &   1.790281e+02   &   4.587219e+02   &   3.934139e+03   &   3.455683e+04   &   8.903309e+05\\
voro$_{5}$   &   2.279882e+01   &   1.422898e+04   &   1.236364e+04   &   3.240000e+05   &   4.041216e+05\\
voro$_{6}$   &   1.186919e+03   &   3.205795e+04   &   1.897652e+04   &   4.764303e+05   &   1.190725e+05\\
voro$_{7}$   &   4.383398e+01   &   4.341273e+05   &   2.040519e+04   &   5.071311e+04   &   2.846018e+05\\
voro$_{8}$   &   1.365488e+03   &   3.725217e+03   &   1.948574e+04   &   4.217183e+04   &   7.119716e+05\\
voro$_{9}$   &   9.643209e+01   &   1.961614e+04   &   2.692660e+04   &   1.469992e+06   &   2.024868e+05\\
voro$_{10}$   &   6.926221e+01   &   5.753447e+03   &   1.388655e+05   &   1.483853e+05   &   6.775208e+04\\
\midrule
\multicolumn{6}{c}{{BoomerAMG}}\\
\midrule
$\text{Mesh} \backslash L$ & {$64$} & {$128$} & {$256$} & {$512$} & {$1024$}\\
\midrule
voro$_{3}$   &   1.952994e+00   &   1.859682e+00   &   1.865789e+00   &   2.079296e+00   &   2.252361e+00\\
voro$_{4}$   &   2.598864e+00   &   2.003719e+00   &   2.315558e+00   &   2.249157e+00   &   4.428457e+00\\
voro$_{5}$   &   2.193265e+00   &   2.168652e+00   &   2.313831e+00   &   2.111331e+00   &   2.161860e+00\\
voro$_{6}$   &   2.025835e+00   &   2.325579e+00   &   2.022580e+00   &   2.255899e+00   &   2.067713e+00\\
voro$_{7}$   &   2.058416e+00   &   2.499529e+00   &   2.476255e+00   &   3.031269e+00   &   3.160480e+00\\
voro$_{8}$   &   2.133230e+00   &   4.140294e+00   &   2.145998e+00   &   3.188217e+00   &   3.254404e+00\\
voro$_{9}$   &   2.380230e+00   &   2.414223e+00   &   2.361048e+00   &   2.573199e+00   &   2.749394e+00\\
voro$_{10}$   &   2.865473e+00   &   2.348662e+00   &   3.849749e+00   &   2.817626e+00   &   3.012861e+00\\
\midrule
\multicolumn{6}{c}{{ML}}\\
\midrule
$\text{Mesh} \backslash L$ & {$64$} & {$128$} & {$256$} & {$512$} & {$1024$}\\
\midrule
voro$_{3}$   &   1.126942e+02   &   1.322346e+06   &   8.304596e+04   &   2.661181e+05   &   4.677648e+05\\
voro$_{4}$   &   6.514475e+02   &   3.695151e+03   &   9.427280e+03   &   4.154313e+05   &   9.906185e+06\\
voro$_{5}$   &   2.112337e+03   &   1.793968e+04   &   7.138982e+04   &   5.119100e+05   &   7.010341e+05\\
voro$_{6}$   &   3.322679e+01   &   2.971126e+05   &   9.321420e+04   &   3.202032e+06   &   7.281148e+05\\
voro$_{7}$   &   4.312113e+02   &   2.113862e+06   &   1.441598e+05   &   2.152576e+05   &   1.790796e+06\\
voro$_{8}$   &   1.940853e+04   &   3.850792e+04   &   4.735274e+04   &   9.441339e+04   &   6.931470e+05\\
voro$_{9}$   &   3.388850e+02   &   6.846527e+04   &   1.116631e+05   &   1.962133e+06   &   9.194373e+05\\
voro$_{10}$   &   1.472389e+03   &   2.001778e+04   &   3.359400e+05   &   5.141470e+05   &   1.535189e+05\\
\bottomrule
\end{tabular}
\end{table}

\subsection{Performance Comparison of Direct and Iterative Solvers}
\label{computational_times_fig}

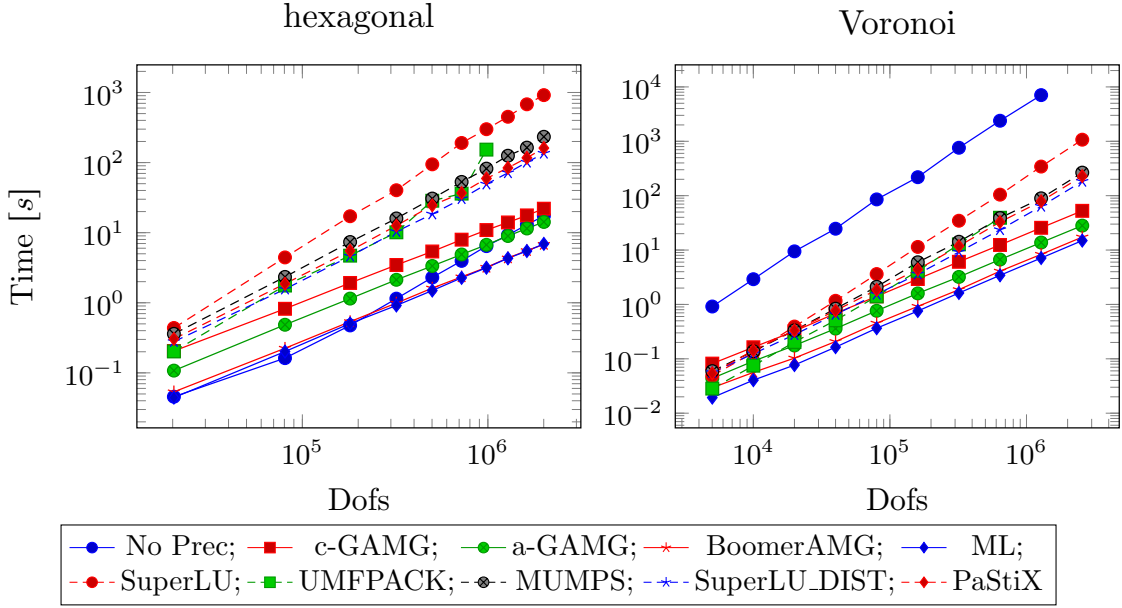
\begin{figure}
	\centering
	\begin{tikzpicture}[scale=1.2][trim axis left]
%
%
	\begin{loglogaxis}[small, xlabel=Dofs, ylabel={Time $[s]$}, title = {hexagonal}, 
	legend columns=5,
	legend entries={$\text{No Prec}$;, $\text{c-GAMG}$;, $\text{a-GAMG}$;, $\text{BoomerAMG}$;, $\text{ML}$;, $\text{SuperLU}$;, $\text{UMFPACK}$;, $\text{MUMPS}$;, $\text{SuperLU\_DIST}$;, $\text{PaStiX}$},
	legend to name=named-hexa-voro,
	cycle list name=mycolorlist]
	\addplot coordinates {
		(20304, 0.045653)
		(80604, 0.16226)
		(180904, 0.47489)
		(321204, 1.1451)
		(501504, 2.2957)
		(721804, 3.9597)
		(982104, 6.4281)
		(1282404, 9.2533)
		(1622704, 12.969)
		(2003004, 17.36)
	};
	\addplot coordinates {
		(20304, 0.20505)
		(80604, 0.82052)
		(180904, 1.9168)
		(321204, 3.458)
		(501504, 5.3899)
		(721804, 7.951)
		(982104, 10.892)
		(1282404, 13.988)
		(1622704, 17.638)
		(2003004, 21.989)
	};
	\addplot coordinates {
		(20304, 0.10766)
		(80604, 0.48551)
		(180904, 1.147)
		(321204, 2.127)
		(501504, 3.3479)
		(721804, 4.8896)
		(982104, 6.7115)
		(1282404, 8.9592)
		(1622704, 11.472)
		(2003004, 14.22)
	};
	\addplot coordinates {
		(20304, 0.053521)
		(80604, 0.22416)
		(180904, 0.53778)
		(321204, 0.9893)
		(501504, 1.619)
		(721804, 2.3374)
		(982104, 3.2013)
		(1282404, 4.3529)
		(1622704, 5.5429)
		(2003004, 6.7944)
	};
	\addplot coordinates {
		(20304, 0.04433)
		(80604, 0.20129)
		(180904, 0.50579)
		(321204, 0.91425)
		(501504, 1.494)
		(721804, 2.266)
		(982104, 3.1606)
		(1282404, 4.3112)
		(1622704, 5.4728)
		(2003004, 6.92)
	};
	\addplot coordinates {
		(20304, 0.43678)
		(80604, 4.4357)
		(180904, 17.16)
		(321204, 40.319)
		(501504, 94.402)
		(721804, 190.54)
		(982104, 299.4)
		(1282404, 450.36)
		(1622704, 678.48)
		(2003004, 917.19)
	};
	\addplot coordinates {
		(20304, 0.20205)
		(80604, 1.76)
		(180904, 4.7001)
		(321204, 10.114)
		(501504, 28.546)
		(721804, 36.074)
		(982104, 153.1)
	};
	\addplot coordinates {
		(20304, 0.36014)
		(80604, 2.3383)
		(180904, 7.3608)
		(321204, 15.961)
		(501504, 30.678)
		(721804, 52.996)
		(982104, 82.006)
		(1282404, 125.26)
		(1622704, 163.5)
		(2003004, 233.31)
	};
	\addplot coordinates {
		(20304, 0.27888)
		(80604, 1.5642)
		(180904, 4.614)
		(321204, 10.364)
		(501504, 18.414)
		(721804, 30.481)
		(982104, 49.179)
		(1282404, 71.369)
		(1622704, 100.93)
		(2003004, 135.5)
	};
	\addplot coordinates {
		(20304, 0.30496)
		(80604, 1.8645)
		(180904, 5.512)
		(321204, 12.598)
		(501504, 24.532)
		(721804, 36.777)
		(982104, 59.314)
		(1282404, 84.039)
		(1622704, 116.91)
		(2003004, 161.21)
	};
	\end{loglogaxis}
	\end{tikzpicture}
	\begin{tikzpicture}[scale=1.2][trim axis left]
	\begin{loglogaxis}[small,
	xlabel=Dofs, 
	title = {Voronoi},
	cycle list name=mycolorlist]
	\addplot coordinates {
		(5006, 0.91552)
		(10008, 2.9097)
		(20007, 9.4653)
		(40011, 24.67)
		(80007, 85.39)
		(160028, 218.74)
		(320020, 760.42)
		(640035, 2387.3)
		(1280053, 7067.2)
	};
	\addplot coordinates {
		(5006, 0.08162)
		(10008, 0.16309)
		(20007, 0.32221)
		(40011, 0.6834)
		(80007, 1.4079)
		(160028, 2.9335)
		(320020, 6.0827)
		(640035, 12.407)
		(1280053, 25.664)
		(2560065, 52.643)
	};
	\addplot coordinates {
		(5006, 0.043763)
		(10008, 0.09054)
		(20007, 0.17557)
		(40011, 0.35923)
		(80007, 0.76393)
		(160028, 1.5958)
		(320020, 3.2037)
		(640035, 6.6865)
		(1280053, 13.743)
		(2560065, 28.153)
	};
	\addplot coordinates {
		(5006, 0.029351)
		(10008, 0.05628)
		(20007, 0.10303)
		(40011, 0.20794)
		(80007, 0.45199)
		(160028, 0.92067)
		(320020, 1.8673)
		(640035, 4.0881)
		(1280053, 8.3346)
		(2560065, 17.311)
	};
	\addplot coordinates {
		(5006, 0.019332)
		(10008, 0.040488)
		(20007, 0.076759)
		(40011, 0.1638)
		(80007, 0.36494)
		(160028, 0.75699)
		(320020, 1.6226)
		(640035, 3.4304)
		(1280053, 7.1346)
		(2560065, 14.924)
	};
	\addplot coordinates {
		(5006, 0.048909)
		(10008, 0.13066)
		(20007, 0.39211)
		(40011, 1.1661)
		(80007, 3.6194)
		(160028, 11.41)
		(320020, 34.604)
		(640035, 104.42)
		(1280053, 342.92)
		(2560065, 1070.3)
	};
	\addplot coordinates {
		(5006, 0.028318)
		(10008, 0.074702)
		(20007, 0.20329)
		(40011, 0.51084)
		(80007, 1.4038)
		(160028, 4.1462)
		(320020, 12.708)
		(640035, 39.266)
	};
	\addplot coordinates {
		(5006, 0.059363)
		(10008, 0.13736)
		(20007, 0.3329)
		(40011, 0.83607)
		(80007, 2.1098)
		(160028, 5.9635)
		(320020, 14.214)
		(640035, 38.819)
		(1280053, 89.477)
		(2560065, 265.41)
	};
	\addplot coordinates {
		(5006, 0.055963)
		(10008, 0.12392)
		(20007, 0.27823)
		(40011, 0.65385)
		(80007, 1.5334)
		(160028, 3.624)
		(320020, 9.1874)
		(640035, 23.621)
		(1280053, 63.586)
		(2560065, 182.87)
	};
	\addplot coordinates {
		(5006, 0.054207)
		(10008, 0.14279)
		(20007, 0.33321)
		(40011, 0.76909)
		(80007, 1.8535)
		(160028, 4.398)
		(320020, 11.807)
		(640035, 32.741)
		(1280053, 78.915)
		(2560065, 228.65)
	};
	\end{loglogaxis}
	\end{tikzpicture}
	\\
	\ref{named-hexa-voro}
	\caption{Time (for setup and solving) 
		of direct and iterative solvers in a loglog scale for the {\em hexagonal meshes} of Table~\ref{tab:mesh-hexa} and for the {\em Voronoi meshes} of Table~\ref{tab:mesh-voro}. 
		With mesh voro$_{10}$, CG without preconditioning requires more than $10000$ iterations to reach convergence.
With meshes voro$_9$, voro$_{10}$ and hexa$_8$, hexa$_9$, hexa$_{10}$ UMFPACK fails due to an \emph{out of memory} error.
	}
	\label{fig:k1-hexa-voro}
\end{figure}

\begin{figure}
	\centering
	\begin{tikzpicture}[scale=1.2][trim axis left]
	\begin{loglogaxis}[
	small,
	xlabel=Dofs, ylabel={Time $[s]$}, title = {horses},
	legend columns=5,
	legend entries={$\text{No Prec}$;, $\text{c-GAMG}$;, $\text{a-GAMG}$;, $\text{BoomerAMG}$;, $\text{ML}$;, $\text{SuperLU}$;, $\text{UMFPACK}$;, $\text{MUMPS}$;, $\text{SuperLU\_DIST}$;, $\text{PaStiX}$},
	legend to name=named-horse-koch, 
	cycle list name=mycolorlist]
	\addplot coordinates {
		(94401, 149.33)
		(135481, 250.89)
		(183961, 380.27)
		(239841, 549.03)
		(303121, 766.32)
		(373801, 1026)
		(451881, 1295.9)
		(537361, 1663)
		(630241, 2081.9)
		(730521, 2571.2)
	};
	\addplot coordinates {
		(94401, 24.673)
		(135481, 38.347)
		(183961, 57.833)
		(239841, 81.981)
		(303121, 116.91)
		(373801, 145.49)
		(451881, 188.48)
		(537361, 242.24)
		(630241, 302.77)
		(730521, 367.86)
	};
	\addplot coordinates {
		(94401, 16.6)
		(135481, 24.273)
		(183961, 33.713)
		(239841, 44.681)
		(303121, 57.57)
		(373801, 71.957)
		(451881, 86.813)
		(537361, 104.55)
		(630241, 123.07)
		(730521, 143.51)
	};
	\addplot coordinates {
		(94401, 11.332)
		(135481, 16.156)
		(183961, 22.014)
		(239841, 30.306)
		(303121, 36.166)
		(373801, 46.732)
		(451881, 57.419)
		(537361, 68.282)
		(630241, 79.241)
		(730521, 91.991)
	};
	\addplot coordinates {
		(94401, 4.5034)
		(135481, 6.6268)
		(183961, 9.3621)
		(239841, 12.751)
		(303121, 16.547)
		(373801, 20.957)
		(451881, 25.622)
		(537361, 30.879)
		(630241, 35.882)
		(730521, 43.889)
	};
	\addplot coordinates {
		(94401, 280.2)
		(135481, 578.07)
		(183961, 1068.7)
		(239841, 1819.5)
		(303121, 2910.7)
		(373801, 4425.4)
		(451881, 7690.3)
		(537361, 8554.2)
	};
	\addplot coordinates {
		(94401, 53.079)
		(135481, 97.411)
		(183961, 151.42)
		(239841, 237.75)
		(303121, 352.39)
	};
	\addplot coordinates {
		(94401, 28.767)
		(135481, 48.24)
		(183961, 75.912)
		(239841, 110.89)
		(303121, 162.01)
		(373801, 141.7)
		(451881, 284.76)
		(537361, 338.76)
		(630241, 452.43)
		(730521, 597)
	};
	\addplot coordinates {
		(94401, 14.244)
		(135481, 27.227)
		(183961, 39.155)
		(239841, 60.67)
		(303121, 78.628)
		(373801, 109.2)
		(451881, 155.55)
		(537361, 198.77)
		(630241, 258.73)
		(730521, 308.16)
	};
	\addplot coordinates {
		(94401, 20.271)
		(135481, 26.959)
		(183961, 44.552)
		(239841, 71.129)
		(303121, 90.443)
		(373801, 125.58)
		(451881, 181.73)
		(537361, 242.53)
		(630241, 284.98)
		(730521, 402.3)
	};
	\end{loglogaxis}
	\end{tikzpicture}
%
%
	\begin{tikzpicture}[scale=1.2][trim axis left]
	\begin{loglogaxis}[
	small,
	xlabel=Dofs, 
	 title ={Koch snowflakes},
	cycle list name=mycolorlist]
	\addplot coordinates {
		(27973, 21.531)
		(49713, 51.627)
		(77661, 97.942)
		(111817, 166.27)
		(152181, 262.51)
		(198753, 327.79)
		(251533, 454.06)
		(310521, 587.12)
		(375717, 761.78)
		(447121, 913.9)
		(524733, 1113.7)
		(608553, 1332.8)
		(698581, 1584.1)
	};
	\addplot coordinates {
		(27973, 16.394)
		(49713, 32.685)
		(77661, 56.046)
		(111817, 86.023)
		(152181, 123.2)
		(198753, 168.41)
		(251533, 226.3)
		(310521, 291.17)
		(375717, 375.47)
		(447121, 443.58)
		(524733, 534.65)
		(608553, 643.22)
		(698581, 745.59)
	};
	\addplot coordinates {
		(27973, 12.55)
		(49713, 23.781)
		(77661, 36.768)
		(111817, 53.713)
		(152181, 73.848)
		(198753, 96.782)
		(251533, 122.69)
		(310521, 153.01)
		(375717, 184.01)
		(447121, 219.29)
		(524733, 258.17)
		(608553, 299.8)
		(698581, 345.2)
	};
	\addplot coordinates {
		(27973, 3.275)
		(49713, 5.5831)
		(77661, 8.7684)
		(111817, 12.42)
		(152181, 16.474)
		(198753, 21.117)
		(251533, 26.733)
		(310521, 32.799)
		(375717, 39.155)
		(447121, 45.871)
		(524733, 53.606)
		(608553, 62.004)
		(698581, 70.977)
	};
	\addplot coordinates {
		(27973, 1.0097)
		(49713, 1.8667)
		(77661, 3.0388)
		(111817, 4.6977)
		(152181, 6.5114)
		(198753, 8.9937)
		(251533, 11.435)
		(310521, 13.811)
		(375717, 17.606)
		(447121, 21.434)
		(524733, 24.59)
		(608553, 29.646)
		(698581, 33.412)
	};
	\addplot coordinates {
		(27973, 1.7345)
		(49713, 3.198)
		(77661, 5.2957)
		(111817, 8.1646)
		(152181, 11.924)
		(198753, 16.872)
		(251533, 23.069)
		(310521, 30.703)
	};
	\addplot coordinates {
		(27973, 1.3959)
		(49713, 2.5131)
		(77661, 3.9955)
		(111817, 5.8966)
		(152181, 8.0952)
		(198753, 10.863)
		(251533, 13.709)
		(310521, 16.57)
		(375717, 20.464)
		(447121, 24.926)
	};
	\addplot coordinates {
		(27973, 0.94683)
		(49713, 1.4538)
		(77661, 2.5893)
		(111817, 3.7623)
		(152181, 5.093)
		(198753, 5.8099)
		(251533, 7.16)
		(310521, 10.372)
		(375717, 12.739)
		(447121, 12.806)
		(524733, 17.894)
		(608553, 20.635)
		(698581, 20.187)
	};
	\addplot coordinates {
		(27973, 1.0734)
		(49713, 1.9793)
		(77661, 3.0891)
		(111817, 4.5162)
		(152181, 6.1561)
		(198753, 7.2941)
		(251533, 9.2413)
		(310521, 12.467)
		(375717, 13.956)
		(447121, 16.505)
		(524733, 19.547)
		(608553, 22.45)
		(698581, 28.09)
	};
	\addplot coordinates {
		(27973, 1.3945)
		(49713, 2.5254)
		(77661, 4.6091)
		(111817, 6.5096)
		(152181, 9.0957)
		(198753, 10.163)
		(251533, 12.863)
		(310521, 15.893)
		(375717, 19.586)
		(447121, 23.142)
		(524733, 27.031)
		(608553, 37.387)
		(698581, 42.312)
	};
	\end{loglogaxis}
	\end{tikzpicture}
	\ref{named-horse-koch}
	\caption{	Time (for setup and solving) 
		of direct and iterative solvers in a loglog scale for \emph{meshes of horses} (Table~\ref{tab:mesh-horse}) and for meshes with {\em Koch snowflakes} (Table~\ref{tab:mesh-koch}). 
		SuperLU fails due to an \emph{out of memory} error on meshes horse$_9$ and horse$_{10}$ and koch$_9$--koch$_{13}$ whereas UMFPACK fails on meshes horse$_6$--horse$_{10}$ 
			 and 
			 koch$_{11}$--koch$_{13}$.
	}
	\label{fig:k1-koch}
\end{figure}
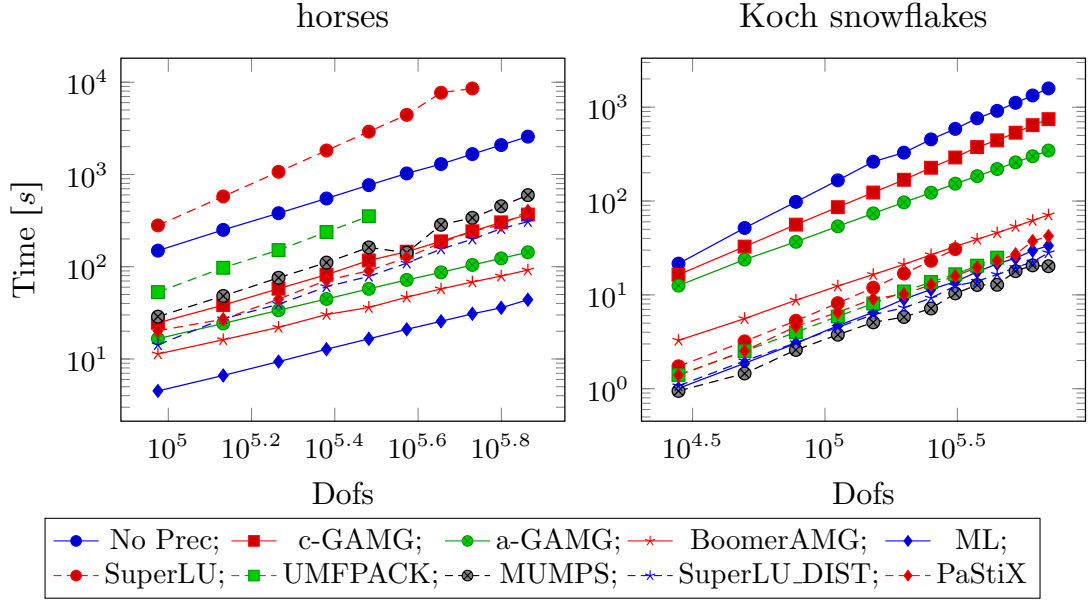

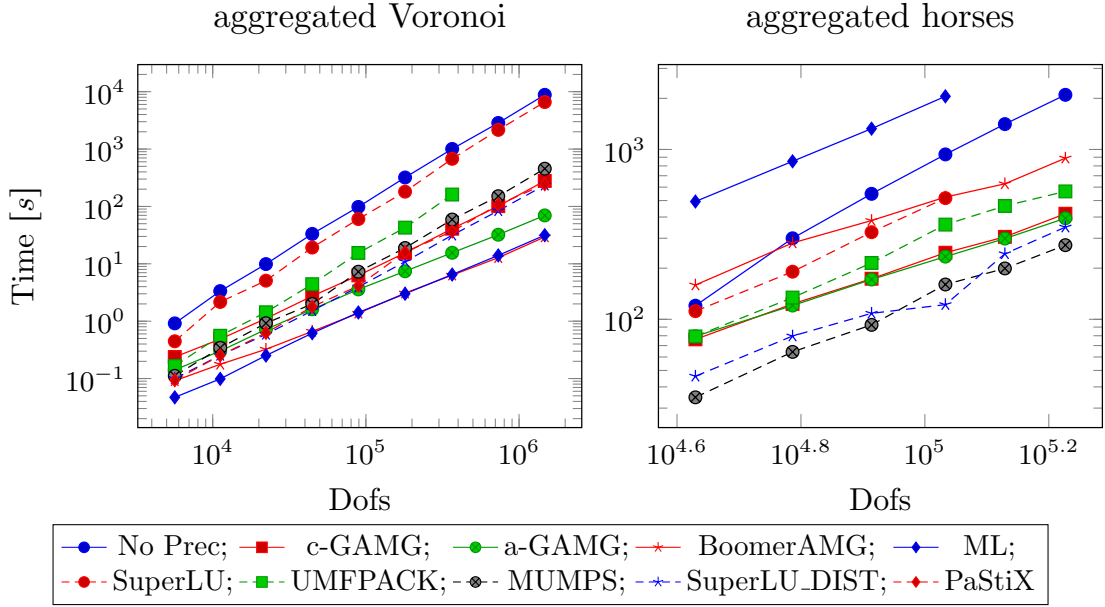
\begin{figure}
	\centering
	\begin{tikzpicture}[scale=1.2][trim axis left]
	\begin{loglogaxis}[
	small,
	xlabel=Dofs, ylabel={Time $[s]$}, title = {aggregated Voronoi}, 
	legend columns=5,
	legend entries={$\text{No Prec}$;, $\text{c-GAMG}$;, $\text{a-GAMG}$;, $\text{BoomerAMG}$;, $\text{ML}$;, $\text{SuperLU}$;, $\text{UMFPACK}$;, $\text{MUMPS}$;, $\text{SuperLU\_DIST}$;, $\text{PaStiX}$},
	legend to name=named-aggr, 
	cycle list name=mycolorlist]
	\addplot coordinates {
		(5652, 0.91591)
		(11190, 3.3472)
		(22250, 9.8986)
		(44615, 33.218)
		(89283, 98.114)
		(180219, 319.52)
		(365921, 1004.9)
		(731552, 2850)
		(1471656, 8816)
	};
	\addplot coordinates {
		(5652, 0.2376)
		(11190, 0.48695)
		(22250, 1.1118)
		(44615, 2.7327)
		(89283, 6.2733)
		(180219, 15.513)
		(365921, 40.934)
		(731552, 101.21)
		(1471656, 278.36)
	};
	\addplot coordinates {
		(5652, 0.14155)
		(11190, 0.30161)
		(22250, 0.7198)
		(44615, 1.611)
		(89283, 3.5809)
		(180219, 7.4151)
		(365921, 15.585)
		(731552, 32.126)
		(1471656, 69.776)
	};
	\addplot coordinates {
		(5652, 0.091514)
		(11190, 0.17726)
		(22250, 0.32378)
		(44615, 0.66466)
		(89283, 1.3737)
		(180219, 3.0861)
		(365921, 6.3442)
		(731552, 12.776)
		(1471656, 29.381)
	};
	\addplot coordinates {
		(5652, 0.047074)
		(11190, 0.098826)
		(22250, 0.25157)
		(44615, 0.61792)
		(89283, 1.4189)
		(180219, 3.0051)
		(365921, 6.5038)
		(731552, 14.012)
		(1471656, 31.551)
	};
	\addplot coordinates {
		(5652, 0.44655)
		(11190, 2.1618)
		(22250, 5.0852)
		(44615, 19.197)
		(89283, 60.273)
		(180219, 180.88)
		(365921, 676.15)
		(731552, 2160.1)
		(1471656, 6601.2)
	};
	\addplot coordinates {
		(5652, 0.16557)
		(11190, 0.55736)
		(22250, 1.431)
		(44615, 4.4137)
		(89283, 15.475)
		(180219, 42.723)
		(365921, 161.19)
	};
	\addplot coordinates {
		(5652, 0.11158)
		(11190, 0.34276)
		(22250, 0.92242)
		(44615, 1.9924)
		(89283, 7.2822)
		(180219, 18.892)
		(365921, 59.134)
		(731552, 150.93)
		(1471656, 454.94)
	};
	\addplot coordinates {
		(5652, 0.09993)
		(11190, 0.25352)
		(22250, 0.58947)
		(44615, 1.521)
		(89283, 3.9478)
		(180219, 10.892)
		(365921, 31.527)
		(731552, 84.079)
		(1471656, 234.54)
	};
	\addplot coordinates {
		(5652, 0.093522)
		(11190, 0.25285)
		(22250, 0.62481)
		(44615, 1.7812)
		(89283, 4.1225)
		(180219, 16.042)
		(365921, 36.667)
		(731552, 107.36)
		(1471656, 243.77)
	};
	\end{loglogaxis}
	\end{tikzpicture}
%
%
%
	\begin{tikzpicture}[scale=1.2][trim axis left]
	\begin{loglogaxis}[
	small,
	xlabel=Dofs, title = {aggregated horses}, 
	cycle list name=mycolorlist]
	\addplot coordinates {
		(42621, 120.26)
		(61151, 299.89)
		(82007, 547.79)
		(107899, 934.99)
		(134589, 1410.2)
		(168525, 2097)
	};
	\addplot coordinates {
		(42621, 76.214)
		(61151, 123.16)
		(82007, 173.28)
		(107899, 246.47)
		(134589, 305.76)
		(168525, 418.23)
	};
	\addplot coordinates {
		(42621, 79.454)
		(61151, 120.5)
		(82007, 171.53)
		(107899, 233.86)
		(134589, 298.65)
		(168525, 393.51)
	};
	\addplot coordinates {
		(42621, 159.1)
		(61151, 281.54)
		(82007, 381.52)
		(107899, 524.13)
		(134589, 627.53)
		(168525, 888.74)
	};
	\addplot coordinates {
		(42621, 493.71)
		(61151, 851.44)
		(82007, 1325.6)
		(107899, 2055)
	};
	\addplot coordinates {
		(42621, 111.62)
		(61151, 190.7)
		(82007, 325.71)
		(107899, 516.99)
	};
	\addplot coordinates {
		(42621, 79.294)
		(61151, 134.25)
		(82007, 214.58)
		(107899, 361.02)
		(134589, 464.82)
		(168525, 566.2)
	};
	\addplot coordinates {
		(42621, 34.721)
		(61151, 64.283)
		(82007, 92.711)
		(107899, 160.29)
		(134589, 199.31)
		(168525, 273.21)
	};
	\addplot coordinates {
		(42621, 46.22)
		(61151, 79.695)
		(82007, 108.41)
		(107899, 121.97)
		(134589, 242.98)
		(168525, 348.91)
	};
	\end{loglogaxis}
	\end{tikzpicture}
	\ref{named-aggr} 
	\caption{
		Time (for setup and solving) 
		of direct and iterative solvers in a loglog scale for meshes of \emph{aggregated Voronoi} and  \emph{aggregated horses} (Tables~\ref{tab:mesh-voro-agg}-\ref{tab:mesh-horse-agg}). 
		Both SuperLU and UMFPACK fails on a-horse$_5$ and a-horse$_6$ due to an \emph{out of memory} error.}
	\label{fig:k1-horse-agg}
\end{figure}

\begin{figure}
	\centering
	\begin{tikzpicture}[scale=0.8][trim axis left]
	\begin{loglogaxis}[
	small, xlabel = Dofs, title = {$L = 64$}, 
	ylabel={Time $[s]$}, 
	legend columns=5,
	legend entries={$\text{c-GAMG}$;, $\text{a-GAMG}$;, $\text{BoomerAMG}$;, $\text{ML}$;, $\text{SuperLU}$;, $\text{UMFPACK}$;, $\text{MUMPS}$;, $\text{SuperLU\_DIST}$;, $\text{PaStiX}$},
	legend to name=robust-irho3-hexa-time,
	cycle list name=mycolorlist]
	\addplot coordinates {
		(80604, 0.91213)
		(321204, 3.9874)
		(721804, 8.3183)
		(1282404, 15.22)
	};
	\addplot coordinates {
		(80604, 0.60452)
		(321204, 2.9554)
		(721804, 6.0131)
		(1282404, 11.719)
	};
	\addplot coordinates {
		(80604, 0.26308)
		(321204, 1.0634)
		(721804, 2.6356)
		(1282404, 4.5659)
		(2003004, 7.1908)
	};
	\addplot coordinates {
		(80604, 0.306)
		(321204, 1.6689)
		(721804, 3.449)
		(1282404, 6.2838)
		(2003004, 12.95)
	};
	\addplot coordinates {
		(80604, 4.3881)
		(321204, 39.8)
		(721804, 188.9)
		(1282404, 445.33)
		(2003004, 907.64)
	};
	\addplot coordinates {
		(80604, 1.7315)
		(321204, 10.098)
		(721804, 35.454)
	};
	\addplot coordinates {
		(80604, 2.3106)
		(321204, 15.643)
		(721804, 52.596)
		(1282404, 124.48)
		(2003004, 231.43)
	};
	\addplot coordinates {
		(80604, 1.6881)
		(321204, 10.841)
		(721804, 30.718)
		(1282404, 72.509)
		(2003004, 138.22)
	};
	\addplot coordinates {
		(80604, 1.8712)
		(321204, 12.499)
		(721804, 36.769)
		(1282404, 84.145)
		(2003004, 161.01)
	};
	\end{loglogaxis}
	\end{tikzpicture}
	\begin{tikzpicture}[scale=0.8][trim axis left]
	\begin{loglogaxis}[
	small,  xlabel = Dofs, title = {$L = 256$},
	cycle list name=mycolorlist]
	\addplot coordinates {
		(80604, 1.1396)
		(321204, 5.5811)
		(721804, 13.857)
		(1282404, 27.385)
		(2003004, 32.243)
	};
	\addplot coordinates {
		(80604, 1.5927)
		(321204, 7.9401)
		(721804, 17.538)
		(1282404, 38.986)
		(2003004, 40.656)
	};
	\addplot coordinates {
		(80604, 0.27117)
		(321204, 1.101)
		(721804, 2.5346)
		(1282404, 4.2508)
		(2003004, 6.8684)
	};
	\addplot coordinates {
		(80604, 1.1405)
		(321204, 5.5657)
		(721804, 12.916)
		(1282404, 30.537)
		(2003004, 36.397)
	};
	\addplot coordinates {
		(80604, 4.3793)
		(321204, 39.85)
		(721804, 189.11)
		(1282404, 445.29)
		(2003004, 907.12)
	};
	\addplot coordinates {
		(80604, 3.1394)
		(321204, 14.702)
		(721804, 36.378)
	};
	\addplot coordinates {
		(80604, 2.3095)
		(321204, 15.746)
		(721804, 52.616)
		(1282404, 124.58)
		(2003004, 232.23)
	};
	\addplot coordinates {
		(80604, 1.6736)
		(321204, 10.779)
		(721804, 31.157)
		(1282404, 72.302)
		(2003004, 137.11)
	};
	\addplot coordinates {
		(80604, 1.8742)
		(321204, 12.279)
		(721804, 36.772)
		(1282404, 84.136)
		(2003004, 161.23)
	};
	\end{loglogaxis}
	\end{tikzpicture}
	\begin{tikzpicture}[scale=0.8][trim axis left]
	\begin{loglogaxis}[
	small,
	 xlabel = Dofs, title = {$L = 1024$},
	cycle list name=mycolorlist]
	\addplot coordinates {
		(80604, 5.0487)
		(321204, 14.182)
		(721804, 56.76)
		(1282404, 65.313)
		(2003004, 136.96)
	};
	\addplot coordinates {
		(80604, 14.688)
		(321204, 33.694)
		(721804, 89.818)
		(1282404, 135.07)
		(2003004, 293.3)
	};
	\addplot coordinates {
		(80604, 0.28678)
		(321204, 1.2879)
		(721804, 2.608)
		(1282404, 4.7267)
		(2003004, 6.9779)
	};
	\addplot coordinates {
		(80604, 16.315)
		(321204, 32.948)
		(721804, 93.934)
		(1282404, 119.62)
		(2003004, 326.04)
	};
	\addplot coordinates {
		(80604, 4.3622)
		(321204, 39.667)
		(721804, 189.55)
		(1282404, 445.16)
		(2003004, 906.21)
	};
	\addplot coordinates {
		(80604, 5.8249)
		(321204, 23.247)
		(721804, 218.18)
	};
	\addplot coordinates {
		(80604, 2.3098)
		(321204, 15.696)
		(721804, 52.661)
		(1282404, 124.6)
		(2003004, 231.67)
	};
	\addplot coordinates {
		(80604, 1.6676)
		(321204, 10.846)
		(721804, 31.021)
		(1282404, 72.557)
		(2003004, 137.36)
	};
	\addplot coordinates {
		(80604, 1.8786)
		(321204, 12.432)
		(721804, 36.804)
		(1282404, 83.913)
		(2003004, 160.99)
	};
	\end{loglogaxis}
	\end{tikzpicture}
	\ref{robust-irho3-hexa-time}
	\caption{
		Time (for setup and solving) 
		of direct and iterative solvers for the hexagonal meshes hexa$_2$, hexa$_4$, hexa$_6$, hexa$_8$ and hexa$_{10}$ of  Table~\ref{tab:mesh-hexa} partitioned 
		in $L = 64, 256, 1024$ parts. Every elements of a part is assigned the same diffusion coefficient $\rho = 10^\alpha$, with $\alpha\in[-5,5]$. 
}
	\label{fig:hexa-random}
\end{figure}
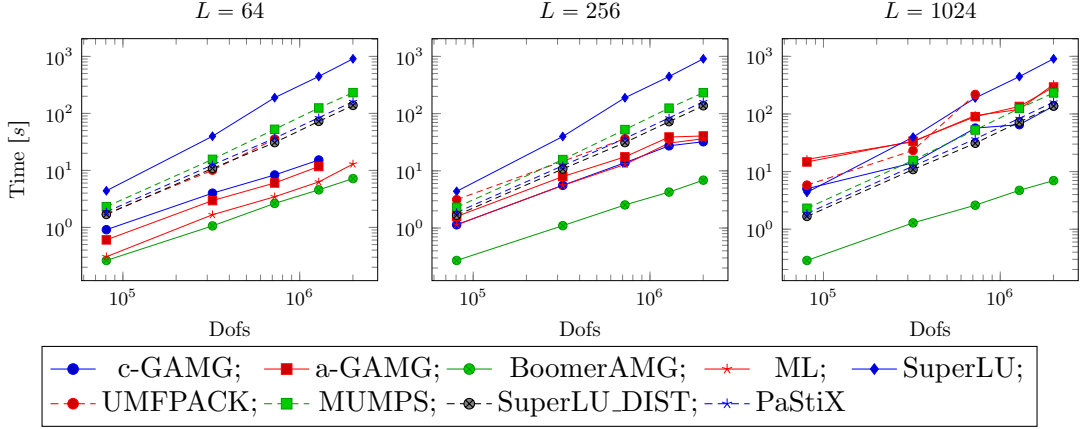

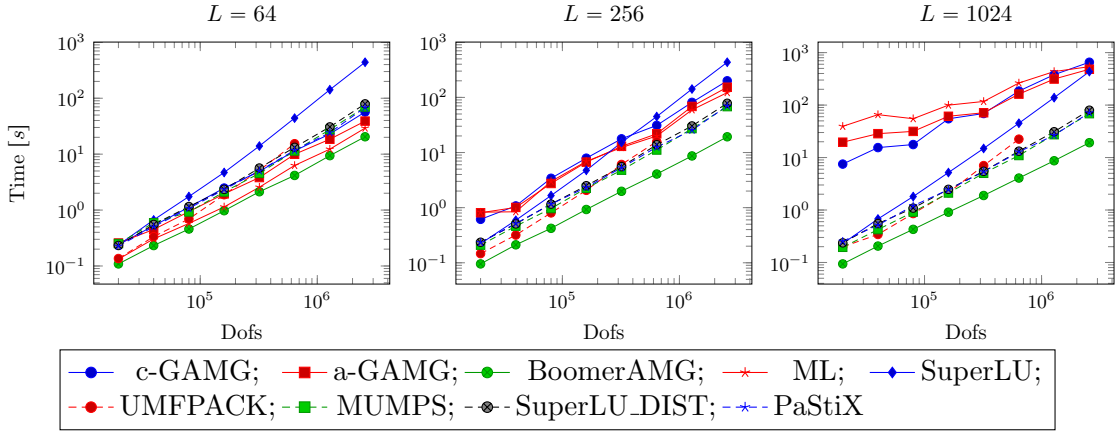
\begin{figure}
	\centering
	\begin{tikzpicture}[scale=0.8][trim axis left]
	\begin{loglogaxis}[
	small,
	xlabel=Dofs, ylabel={Time $[s]$}, title= {$L = 64$},
	legend columns=5,
	legend entries={$\text{c-GAMG}$;, $\text{a-GAMG}$;, $\text{BoomerAMG}$;, $\text{ML}$;, $\text{SuperLU}$;, $\text{UMFPACK}$;, $\text{MUMPS}$;, $\text{SuperLU\_DIST}$;, $\text{PaStiX}$},
	legend to name=robust-irho3-voro-time,
	cycle list name=mycolorlist]
	\addplot coordinates {
		(20007, 0.23919)
		(40011, 0.52486)
		(80007, 1.1242)
		(160028, 2.4742)
		(320020, 4.9139)
		(640035, 11.446)
		(1280053, 22.953)
		(2560065, 56.681)
	};
	\addplot coordinates {
		(20007, 0.25732)
		(40011, 0.44939)
		(80007, 0.93814)
		(160028, 1.948)
		(320020, 3.8165)
		(640035, 9.9439)
		(1280053, 18.485)
		(2560065, 38.991)
	};
	\addplot coordinates {
		(20007, 0.10895)
		(40011, 0.23189)
		(80007, 0.45836)
		(160028, 0.97389)
		(320020, 2.1262)
		(640035, 4.1607)
		(1280053, 9.4007)
		(2560065, 20.469)
	};
	\addplot coordinates {
		(20007, 0.13239)
		(40011, 0.30533)
		(80007, 0.5768)
		(160028, 1.1381)
		(320020, 2.5677)
		(640035, 6.2731)
		(1280053, 12.282)
		(2560065, 29.222)
	};
	\addplot coordinates {
		(20007, 0.23658)
		(40011, 0.6672)
		(80007, 1.7456)
		(160028, 4.7006)
		(320020, 13.965)
		(640035, 44.107)
		(1280053, 141.96)
		(2560065, 440.35)
	};
	\addplot coordinates {
		(20007, 0.13599)
		(40011, 0.3344)
		(80007, 0.69415)
		(160028, 1.8806)
		(320020, 5.149)
		(640035, 15.385)
	};
	\addplot coordinates {
		(20007, 0.25138)
		(40011, 0.5939)
		(80007, 0.9249)
		(160028, 2.114)
		(320020, 4.5343)
		(640035, 11.396)
		(1280053, 27.486)
		(2560065, 69.928)
	};
	\addplot coordinates {
		(20007, 0.23264)
		(40011, 0.5522)
		(80007, 1.1542)
		(160028, 2.3426)
		(320020, 5.6466)
		(640035, 13.214)
		(1280053, 30.713)
		(2560065, 79.336)
	};
	\addplot coordinates {
		(20007, 0.22645)
		(40011, 0.50317)
		(80007, 1.0735)
		(160028, 2.3281)
		(320020, 5.0226)
		(640035, 11.738)
		(1280053, 24.234)
		(2560065, 68.286)
	};
	\end{loglogaxis}
	\end{tikzpicture}
	\begin{tikzpicture}[scale=0.8][trim axis left]
	\begin{loglogaxis}[
	small,
	xlabel=Dofs, title={$L = 256$}, 
	cycle list name=mycolorlist]
	\addplot coordinates {
		(20007, 0.61384)
		(40011, 1.0748)
		(80007, 3.4119)
		(160028, 7.899)
		(320020, 17.728)
		(640035, 31.062)
		(1280053, 81.295)
		(2560065, 200.05)
	};
	\addplot coordinates {
		(20007, 0.80825)
		(40011, 1.0095)
		(80007, 2.7438)
		(160028, 6.6927)
		(320020, 12.949)
		(640035, 21.867)
		(1280053, 68.561)
		(2560065, 153.04)
	};
	\addplot coordinates {
		(20007, 0.095413)
		(40011, 0.21166)
		(80007, 0.42168)
		(160028, 0.92577)
		(320020, 1.9803)
		(640035, 4.0734)
		(1280053, 8.6757)
		(2560065, 19.306)
	};
	\addplot coordinates {
		(20007, 0.82717)
		(40011, 0.84182)
		(80007, 2.9882)
		(160028, 6.9114)
		(320020, 12.307)
		(640035, 19.958)
		(1280053, 59.375)
		(2560065, 122.13)
	};
	\addplot coordinates {
		(20007, 0.22203)
		(40011, 0.5914)
		(80007, 1.6566)
		(160028, 4.7524)
		(320020, 15.251)
		(640035, 44.929)
		(1280053, 141.91)
		(2560065, 434.59)
	};
	\addplot coordinates {
		(20007, 0.1465)
		(40011, 0.31841)
		(80007, 0.80297)
		(160028, 2.0617)
		(320020, 6.0724)
		(640035, 15.798)
	};
	\addplot coordinates {
		(20007, 0.20897)
		(40011, 0.46142)
		(80007, 0.96547)
		(160028, 2.2502)
		(320020, 4.7861)
		(640035, 11.106)
		(1280053, 27.182)
		(2560065, 68.077)
	};
	\addplot coordinates {
		(20007, 0.2357)
		(40011, 0.50776)
		(80007, 1.16)
		(160028, 2.4988)
		(320020, 5.5394)
		(640035, 13.835)
		(1280053, 30.348)
		(2560065, 78.262)
	};
	\addplot coordinates {
		(20007, 0.23317)
		(40011, 0.51628)
		(80007, 1.1609)
		(160028, 2.3292)
		(320020, 5.2971)
		(640035, 12.573)
		(1280053, 25.432)
		(2560065, 70.057)
	};
	\end{loglogaxis}
	\end{tikzpicture}
	\begin{tikzpicture}[scale=0.8][trim axis left]
	\begin{loglogaxis}[
	small,
	xlabel=Dofs, title={$L = 1024$}, 
	cycle list name=mycolorlist]
	\addplot coordinates {
		(20007, 7.5051)
		(40011, 15.481)
		(80007, 17.681)
		(160028, 54.848)
		(320020, 69.004)
		(640035, 185)
		(1280053, 379.29)
		(2560065, 655.66)
	};
	\addplot coordinates {
		(20007, 19.614)
		(40011, 28.51)
		(80007, 31.527)
		(160028, 61.183)
		(320020, 71.455)
		(640035, 161.8)
		(1280053, 313.6)
		(2560065, 481.54)
	};
	\addplot coordinates {
		(20007, 0.093928)
		(40011, 0.20548)
		(80007, 0.42791)
		(160028, 0.90888)
		(320020, 1.8899)
		(640035, 4.0643)
		(1280053, 8.6875)
		(2560065, 19.199)
	};
	\addplot coordinates {
		(20007, 39.306)
		(40011, 65.728)
		(80007, 54.97)
		(160028, 99.723)
		(320020, 117.25)
		(640035, 263.58)
		(1280053, 435.63)
		(2560065, 539.27)
	};
	\addplot coordinates {
		(20007, 0.22929)
		(40011, 0.67522)
		(80007, 1.8051)
		(160028, 5.1732)
		(320020, 14.902)
		(640035, 45.08)
		(1280053, 138.42)
		(2560065, 434.54)
	};
	\addplot coordinates {
		(20007, 0.19838)
		(40011, 0.34193)
		(80007, 0.85573)
		(160028, 2.1392)
		(320020, 7.0265)
		(640035, 22.43)
	};
	\addplot coordinates {
		(20007, 0.19466)
		(40011, 0.42387)
		(80007, 0.91368)
		(160028, 2.1111)
		(320020, 5.0092)
		(640035, 10.934)
		(1280053, 27.262)
		(2560065, 68.618)
	};
	\addplot coordinates {
		(20007, 0.2367)
		(40011, 0.55088)
		(80007, 1.0881)
		(160028, 2.4783)
		(320020, 5.4639)
		(640035, 13.229)
		(1280053, 31.028)
		(2560065, 79.845)
	};
	\addplot coordinates {
		(20007, 0.25418)
		(40011, 0.51623)
		(80007, 1.2039)
		(160028, 2.4435)
		(320020, 5.3808)
		(640035, 12.831)
		(1280053, 26.664)
		(2560065, 71.502)
	};
	\end{loglogaxis}
	\end{tikzpicture}
	\ref{robust-irho3-voro-time}
	\caption{
			Time (for setup and solving) 
		of direct and iterative solvers for the Voronoi meshes voro$_3$--voro$_{10}$ of  Table~\ref{tab:k1-voro}  partitioned 
		in $L = 64, 256, 1024$ parts. Every elements of a part is assigned the same diffusion coefficient $\rho = 10^\alpha$, with $\alpha\in[-5,5]$. 
		}
	\label{fig:voro-random-metis}
\end{figure}

Finally we compare the performance of several direct and iterative solvers available in PETSCs by measuring the elapsed time for both setup and solving.
We recall that, for iterative methods, convergence was decided by the decrease of either the relative or the absolute residual norm as obtained with the direct solver SuperLU\_DIST (see \texttt{rtol} and \texttt{abstol} reported in each Table).
Figures~\ref{fig:k1-hexa-voro},\ref{fig:k1-koch}, \ref{fig:k1-horse-agg},
\ref{fig:hexa-random} and \ref{fig:voro-random-metis} display the time in seconds (for setup and solving) for both direct and direct solvers in a loglog scale. 

The best results are obtained with BoomerAMG and/or with ML for hexagonal, Voronoi and horse meshes, as well as for meshes of aggregated Voronoi cells.
Conversely, when considering Koch snowflake meshes and aggregates of horse cells, 
direct methods outperform AMG preconditioners. 
We recall that 
AMG preconditioners require a high number of iterations for these meshes, as shown in Figure~\ref{fig:k1-koch}-right and  \ref{fig:k1-horse-agg}-right and Tables~\ref{tab:k1-koch}, \ref{tab:k1-horse-agg}, thereby worsening their performance.

Clearly, a more accurate evaluation of the efficiency of AMG techniques for this type of meshes 
 would require a deeper analysis (outside the scope of this work), 
such as an ad-hoc tuning of the AMG parameters or a more accurate analysis of the matrices involved.
%


\section{Conclusion}\label{sec:conclusion}
We numerically investigated the performance of AMG preconditioners
for the solution of a model elliptic problem  on polygonal meshes employing the virtual element method. 
The tested meshes range from the most regular hexagonal to more complex and challenging grids that may mimic features of  more realistic problems such as those based on periodic cellular structures. %
Our results show that CG accelerated multigrid (AMG/CG) is very effective when dealing with either regular (hexagonal) or Voronoi and aggregates of Voronoi meshes.
Moreover, our tests showed that BoomerAMG is also robust if highly varying diffusion coefficients are considered for both set of meshes considered.
Hinging on the results of the present paper, the adoption of AMG preconditioners (in particular BoomerAMG) seems to be a promising approach for solving large linear systems of equations associated with a VEM discretization, in terms of both scalability and reduction of the overall computational cost.

 However, we also verified that, when more complex and challenging meshes are taken into account, such as meshes with many tiny (when compared to the diameter of the elements) edges, most of the considered AMG preconditioners do not preserve scalability, whilst those that retain it, like BoomerAMG or ML, lose most of their efficiency.
 Certainly, the linear system associated with the VEM discretization based on these meshes deserves further investigations that will be addressed in future works.

\bibliographystyle{amsplain}
%
\bibliography{biblio,biblioAMG}

	\end{document}